\documentclass[a4paper,11pt]{article}

\usepackage{amsmath}
\usepackage{amssymb}
\usepackage[applemac]{inputenc}
\usepackage{latexsym}
\usepackage{graphicx}
\usepackage{color}
\usepackage[english]{babel}

\usepackage{mathrsfs}

\DeclareFontFamily{U}{calligra}{}
\DeclareFontShape{U}{calligra}{m}{n}{<->callig15}{}

\newcommand{\calE}{{\!\!\text{\usefont{U}{calligra}{m}{n}E}\,\,}}
\newcommand{\calF}{{\!\!\text{\usefont{U}{calligra}{m}{n}F}\,\,}}

\newcommand{\id}{{\rm id}}
\newcommand{\ad}{{\rm ad}}

\newcommand{\Tr}{{\rm Tr}\,}
\newcommand{\del}{\partial}

\newtheorem{thm}{Theorem}[section]
\newtheorem{prop}[thm]{Proposition}
\newtheorem{rem}[thm]{Remark}
\newtheorem{lem}[thm]{Lemma}
\newtheorem{cor}[thm]{Corollary}

\newtheorem{defe}[thm]{Definition}

\def\proofp{{\hspace{-0.5 cm} $\Box$ \textbf{Proof of Proposition }}}

\def\prooft{{\hspace{-0.5 cm} $\blacksquare$
\textbf{Proof of Theorem }}}
\def\proofl{{\hspace{-0.5 cm} $\vartriangle$ \textbf{Proof of Lemma  }}}
\def\dem{{\hspace{-0.5 cm} $\Box$ \textbf{Proof of Proposition }}}
\def\demcor{{\hspace{-0.5 cm} $\Box$ \textbf{Proof of Corollary  }}}
\def\demthm{{\hspace{-0.5 cm} $\blacksquare$
\textbf{Proof of Theorem  }}}
\def\demlem{{\hspace{-0.5 cm} $\vartriangle$ \textbf{Proof of Lemma }}}

\def\cqfd{{\hfill $\Box$}}
\def\cqfdt{{\hfill $\blacksquare$}}
\def\cq{{\hfill $\vartriangle$}}

\def\R{{\mathbb{R}}}

\def\N{{\mathbb{N}}}

\def\g{{\mathfrak{g}}}

\def\a{{\mathfrak{a}}}
\def\b{{\mathfrak{b}}}

\begin{document}
\pagestyle{plain}
\title{ Totally geodesic submanifolds in the manifold SPD \\of symmetric positive-definite real matrices }

\author{Alice Barbara Tumpach\footnote{{\tt alice-barbora.tumpach@univ-lille.fr}, Institut CNRS Pauli, Vienna, Austria.},\; Gabriel Larotonda\footnote{{\tt glaroton@dm.uba.ar}, UBA \& CONICET, Buenos Aires, Argentina}}

\date{}

\maketitle

\abstract This paper is a self-contained exposition of the geometry of symmetric positive-definite real $n\times n$ matrices $\operatorname{SPD}(n)$, including necessary and sufficent conditions for a submanifold $\mathcal{N} \subset\operatorname{SPD}(n)$ to be totally geodesic for the affine-invariant  Riemannian metric. A non-linear projection $x\mapsto \pi(x)$ on a totally geodesic submanifold is defined. This projection has the minimizing property with respect to the Riemannian metric: it maps an arbitrary point $x \in\operatorname{SPD}(n)$ to the unique closest element $\pi(x)$ in the totally geodesic submanifold  for the distance defined by the  affine-invariant  Riemannian metric. Decompositions of the space $\operatorname{SPD}(n)$ follow, as well as variants of the polar decomposition of non-singular matrices known as Mostow's decompositions. Applications to decompositions of covariant matrices are mentioned.

\noindent{\it Subject Classification:} Riemannian Geometry.

 \tableofcontents

\section{Introduction}
Mostow's decomposition theorem is a generalization of polar decomposition and goes back to 1955 with the original proof by Mostow (\cite{Mostow}). Different proofs were given in \cite{pr,Lar,Conde, Tum} in the context of complex infinite-dimensional Lie groups, in \cite{AL} in the context of von Neumann algebras, and \cite{Miglioli} in the context of Finsler Lie groups. In this paper, we consider the real group $\operatorname{GL}(n, \mathbb{R})$, which is not semi-simple, but appears  in a lot of applications. The aim of this paper is to provide a self-contained exposition of the geometry needed to understand this Theorem, and advertise some possible applications, which should be important in the fields of Probablity and Statistics. For a general exposition of the geometry of the manifold of symmetric positive-definite real matrices we refer to \cite{Bhatia}.

\subsection{Statement of  Mostow's Theorem}

Let us first state the theorem.

\begin{thm}\label{debut}
Let $\operatorname{E}$ be a real subspace of the vector space $\operatorname{Sym}(n)$ of symmetric $n\times n$ real matrices, and denote by $\calE := \exp(\operatorname{E})$ the image of $\operatorname{E}$ by the exponential of matrices
\begin{equation}
\calE = \{y = \exp u :=1+ u + \frac{u^2}{2} + \frac{u^3}{3!} + \cdots, u \in \operatorname{E}\}.
\end{equation}
We endow the manifold $\operatorname{SPD}(n)$ of symmetric positive-definite real $n\times n$ matrices with the Riemannian metric (sometimes called affined Riemannian metric) given at $x\in \operatorname{SPD}(n)$ by:
\begin{equation}\label{def_Riem}
\langle X, Y\rangle_x := \Tr\left(x^{-1}X x^{-1} Y\right) = \Tr \left(x^{-\frac{1}{2}}X x^{-\frac{1}{2}}\cdot x^{-\frac{1}{2}}Y x^{-\frac{1}{2}}\right),
\end{equation} 
where the dot denotes the product of matrices.
\begin{enumerate}
\item
Then $\calE := \exp(\operatorname{E})$ is a totally geodesic submanifold of $\operatorname{SPD}(n)$ and  geodesically convex in $\operatorname{SPD}(n)$,  if and only if the vector space $\operatorname{E}\subset \operatorname{Sym}(n)$ satisfies
\begin{equation}\label{totgeo}
\left[ X, [X, Y]\right] \in \operatorname{E}, \qquad for~all~~~ X, Y \in \operatorname{E},
\end{equation}
where for any matrices $A$, $B$, the bracket is defined as $[A, B] := A\cdot B- B\cdot A$.
\item In this case, any $x\in \operatorname{SPD}(n)$ can be decomposed uniquely as a product $e\cdot f\cdot e$, where 
$e\in \calE := \exp(\operatorname{E})$ and $f\in \calF := \exp(\operatorname{E}^\perp)$.  The matrix $x$ admits a  unique projection $\pi(x)$ onto $\calE := \exp(\operatorname{E})$ which minimises the geodesic distance from $x$ to $\calE := \exp(\operatorname{E})$, i.e. $$\operatorname{dist}(x, \calE) = \operatorname{dist}(x, \pi(x)).$$  The matrix $e$ is related to $\pi(x)$ by
\begin{equation}
e = \pi(x)^{\frac{1}{2}}.
\end{equation}
\item 
If \eqref{totgeo} is satisfied, any $g\in\operatorname{GL}(n, \mathbb{R})$ can be decomposed uniquely as a product $k\cdot \exp V \cdot \exp W$, where $k$ is an orthogonal matrix, $W$ is in $\operatorname{E}$ and $V$ is in the orthogonal space $F := \operatorname{E}^\perp$. The Lie group $\operatorname{GL}(n, \mathbb{R})$ is
diffeomorphic to the product  
\begin{equation}\label{Mostow_dec}
\operatorname{GL}(n, \mathbb{R}) = \operatorname{O}(n) \times \exp \operatorname{F} \times \exp \operatorname{E} = \operatorname{O}(n) \times \calF \times \calE,
\end{equation}
where 
\begin{itemize}
\item$\operatorname{O}(n)$ is the orthogonal group 
\begin{equation}
\operatorname{O}(n) := \{ M \in \operatorname{GL}(n, \mathbb{R}), M^T M = M M^T = \operatorname{Id}\}
\end{equation}
(here $\operatorname{Id}$ is the $n\times n$ identity matrix).
\item $\operatorname{F}$  denotes the orthogonal of $\operatorname{E}$ in $\operatorname{Sym}(n)$, 
\begin{equation}
\operatorname{F} := \operatorname{E}^\perp = \{ Y\in \operatorname{Sym}(n), \Tr X Y = 0~~~\forall X \in \operatorname{E}\},
\end{equation}
\item $\calF := \exp \operatorname{F}.$
\end{itemize}
\end{enumerate}
\end{thm}

\begin{cor}\label{cor_convex}
Any geodesically complete convex submanifold $\mathcal{N}$ of $\operatorname{SPD}(n)$ is of the form $$\mathcal{N} = x \cdot \exp \operatorname E \cdot \,x$$ for some $x\in \operatorname{SPD}(n)$, and some subspace $E \subset \operatorname{Sym}(n)$ satisfying $[X, [X, Y]] \in \operatorname E$ for all $X, Y \in \operatorname E$.
\end{cor}

\subsection{Examples of Applications}

\subsubsection{Geodesic projection on the space of diagonal matrices}
Let $\bold{X} := (X_1, X_2, \cdots, X_n)$ be a real random vector of dimension $n$ and consider its \textit{covariance matrix}
\begin{equation}\label{def_covariance}
\operatorname{cov}(\bold{X}, \bold{X}) = \mathbb{E}\left[\left(\bold{X}-\mathbb{E}(\bold{X})\right)\left(\bold{X}-\mathbb{E}(\bold{X})\right)^T\right]
\end{equation}
where the operator $\mathbb{E}$  denotes the expected value (mean) of its argument.

Let $\operatorname{E}_1\subset \operatorname{Sym}(n)$ be the space of real diagonal matrices. Then condition~\eqref{totgeo} is satisfied, hence by Theorem~\ref{debut} (see Theorem~\ref{projection} for more details), any symmetric positive-definite matrix $x\in \operatorname{SPD}(n)$ has a unique projection $\pi(x)$ in $\calE_1 := \exp(\operatorname{E}_1)$ that minimises (for the Riemannian metric~\eqref{def_Riem}) the distance between $x$ and the manifold of real positive-definite diagonal matrices $\calE_1$. Moreover
\begin{equation}
x = \pi(x)^{\frac{1}{2}}\cdot f \cdot \pi(x)^{\frac{1}{2}},
\end{equation}
for a unique element $f$ in $\,\calF_1:= \exp(\operatorname{E}_1^{\perp})$ where 
$$
\operatorname{E}_1^{\perp} := \left\{ Y = (Y_{ij}) \in \operatorname{Sym}(n), Y_{ii} = 0 \right\}.
$$

\begin{rem}{\rm
The \textit{correlation matrix} of $\bold{X}$ is defined as
\begin{equation}\label{def_correlation}
\operatorname{corr}(\bold{X}, \bold{X}) = \left(\operatorname{diag}(\operatorname{cov}(\bold{X},\bold{X})\right)^{-\frac{1}{2}}\cdot \operatorname{cov}(\bold{X}, \bold{X}) \cdot \left(\operatorname{diag}(\operatorname{cov}(\bold{X},\bold{X})\right)^{-\frac{1}{2}},
\end{equation}
where $\operatorname{diag}(\operatorname{cov}(\bold{X},\bold{X}))$ is the matrix of diagonal elements of $\operatorname{cov}(\bold{X},\bold{X})$, i.e. the diagonal $n\times n$ matrix containing the variances of $X_i$, $i = 1, \dots, n$. 
Then
$$
\operatorname{cov}(\bold{X}, \bold{X}) = \left(\operatorname{diag}(\operatorname{cov}(\bold{X},\bold{X})\right)^{\frac{1}{2}}\cdot \operatorname{corr}(\bold{X}, \bold{X}) \cdot \left(\operatorname{diag}(\operatorname{cov}(\bold{X},\bold{X})\right)^{\frac{1}{2}}.
$$
In this decomposition, we remark that $\operatorname{cov}(\bold{X}, \bold{X})$ is the image of $\operatorname{corr}(\bold{X}, \bold{X})$ by the isometry of $\operatorname{SPD}(n)$ which sends $y$ to 
$$
\left(\operatorname{diag}(\operatorname{cov}(\bold{X},\bold{X})\right)^{\frac{1}{2}}\cdot y \cdot \left(\operatorname{diag}(\operatorname{cov}(\bold{X},\bold{X})\right)^{\frac{1}{2}}.
$$
However, in general, $\operatorname{diag}(x)$ differs from the projection $\pi(x)$ onto the manifold of real positive-definite diagonal matrices for the Riemannian metric~\eqref{def_Riem}. In fact, $\operatorname{diag}(\operatorname{cov}(\bold{X},\bold{X}))$ is the projection of $\operatorname{cov}(\bold{X},\bold{X})$ onto the vector space $\operatorname{E}_1\subset \operatorname{Sym}(n)$ of diagonal matrices for the Euclidian scalar product given by  
\begin{equation}
\begin{array}{llll}
\left(\cdot, \cdot \right): &\operatorname{Sym}(n) \times \operatorname{Sym}(n) & \longrightarrow & \mathbb{R}\\
& (X, Y) & \longmapsto & \Tr(XY).
\end{array}
\end{equation}
}
\end{rem}

\begin{lem} 
Let $x=\pi(x)^{1/2}e^v\pi(x)^{1/2}$ be Mostow's decompositin for $x$, with $v^*=v\in E_1^{\perp}$ a matrix with zeros on the diagonal. Then $\operatorname{diag}(x)=\pi(x)$ if and only if $\operatorname{diag}(e^v)=1$.
\end{lem}
\textbf{Proof}: Write $x=Pe^vP$ for short, and note that the assumption is $P^2=\operatorname{diag}(x)$. Then $P^2=\operatorname{diag}(Pe^vP)$, thus $1=P^{-1}\operatorname{diag}(P e^v P)P^{-1}$. But taking diagonal has the bi-module property, i.e. $\operatorname{diag}(d_1ad_2)=d_1\operatorname{diag}(a)d_2$ for any matrix $a$ and any diagonal matrices $d_1,d_2$. Hence $1=\operatorname{diag}(e^v)$. 

\cqfd

In particular for $n=2$ note that $v$ must be of the form $v=\left(\begin{array}{cc}0 & a \\ a& 0\end{array}\right)$, thus a straighforward computation shows that 
$$
e^v=\left(\begin{array}{cc}\cosh(a) & \sinh(a) \\ \sinh(a)& \cosh(a)\end{array}\right),
$$
and $\operatorname{diag}(e^v)=1$ is only possible if $\cosh(a)=1$, or equivalently if $a=0$. Thus it must be $v=0$ thus  $x=\pi(x)$ was already a diagonal matrix.

\subsubsection{Geodesic projection on the space of block diagonal matrices and on the space of block-anti-diagonal matrices}

Consider now two different random vectors $\bold{X}$ of dimension $p$ and $\bold{Y}$ of dimension $p$, and the random vector $$\bold{Z} := \left(\begin{array}{c} \bold{X}\\ \bold{Y}\end{array}\right)$$ of dimension $p+q =: n$. The \textit{joint covariance matrix}  $\Sigma$ of $\bold{X}$ and $\bold{Y}$ is defined as
\begin{equation}\label{def_joint_covariance}
\Sigma = \mathbb{E}\left[\left(\bold{Z}-\mathbb{E}(\bold{Z})\right)\left(\bold{Z}-\mathbb{E}(\bold{Z})\right)^T\right] = \left(\begin{array}{cc} \operatorname{cov}(\bold{X},\bold{X}) & \operatorname{cov}(\bold{X},\bold{Y})\\ \operatorname{cov}(\bold{X},\bold{Y}) & \operatorname{cov}(\bold{Y},\bold{Y})\end{array}\right),
\end{equation}
where $\operatorname{cov}(\bold{X},\bold{Y}) $ denotes the \textit{cross-covariance matrix} between $\bold{X}$ and $\bold{Y}$:
\begin{equation}\label{def_cross_covariance}
\operatorname{cov}(\bold{X}, \bold{Y}) := \mathbb{E}\left[\left(\bold{X}-\mathbb{E}(\bold{X})\right)\left(\bold{Y}-\mathbb{E}(\bold{Y})\right)^T\right].
\end{equation}

Let $\operatorname{E}_2\subset \operatorname{Sym}(n)$ be the vector space of block-diagonal matrices and 
$\operatorname{E}_3\subset \operatorname{Sym}(n)$ be the vector space of block-anti-diagonal matrices:
$$
\begin{array}{l}
\operatorname{E}_2 = \left\{ Z \in \operatorname{Sym}(n), Z = \left(\begin{array}{cc} X & 0\\ 0 & Y \end{array}\right), X \in \operatorname{Sym}(p), Y \in \operatorname{Sym}(q)\right\}
\\
\operatorname{E}_3 = \left\{ Z \in \operatorname{Sym}(n), Z = \left(\begin{array}{cc} 0 & C\\ C^T & 0 \end{array}\right), C \in  \operatorname{Mat}(p, q)\right\}.
\end{array}
$$
Then $\operatorname{E}_2$ and $\operatorname{E}_3$ both satisfy condition~\eqref{totgeo} and are orthogonal to each other for the Euclidian scalar product on $\operatorname{Sym}(n)$.  Note that $\operatorname{E}_2$ is a Lie algebra, whereas $\operatorname{E}_3$ is not. Therefore the cross-covariance matrix $\Sigma$ has a unique expression of the form
\begin{equation}\label{DAD}
\Sigma = \exp D \cdot \exp A \cdot \exp D,
\end{equation}
where $D$ is block-diagonal and $A$ is block-anti-diagonal, as well as a unique expression of the form
\begin{equation}\label{ADA}
\Sigma = \exp A' \cdot \exp D' \cdot \exp A',
\end{equation}
where $D'$ is block-diagonal and $A'$ is block-anti-diagonal.

\smallskip

Moreover since the even powers of a block-anti-diagonal matrix are diagonal, one has the decomposition
$$
\exp A = \cosh A + \sinh A,
$$
where 
\begin{equation}
\begin{array}{l}
\cosh A := \frac{\left(\exp(A) + \exp(-A)\right)}{2} = \sum_{i = 0}^n \frac{{A}^{2n}}{(2n)!} \in \operatorname{E}_2\\
\sinh A := \frac{\left(\exp(A) - \exp(-A)\right)}{2} = \sum_{i = 1}^n \frac{{A}^{(2n-1)}}{(2n-1)!} \in \operatorname{E}_3\\
\end{array}
\end{equation}
Therefore equation~\eqref{DAD} leads to a decomposition
\begin{equation}
\Sigma = \exp D \cosh A \exp D + \exp D \sinh A \exp D,
\end{equation}
where $\exp D \cosh A \exp D \in   \operatorname{E}_2$ is block-diagonal and $\exp D \sinh A \exp D \in \operatorname{E}_3$ is block-anti-diagonal.
Similarly equation~\eqref{ADA} leads to a decomposition of $\Sigma$ in a sum
\begin{equation}
\left(\cosh A' \cdot \exp D' \cdot \cosh A' + \sinh A' \cdot \exp D' \sinh A'\right) + \left(\sinh A' \exp D' \cosh A' + \cosh A' \exp D' \sinh A' \right)
\end{equation}
where $$\left(\cosh A' \cdot \exp D' \cdot \cosh A' + \sinh A' \cdot \exp D' \sinh A'\right) \in   \operatorname{E}_2$$ is block-diagonal and $$\left(\sinh A' \exp D' \cosh A' + \cosh A' \exp D' \sinh A' \right)\in  \operatorname{E}_3$$ is block-anti-diagonal.

\begin{rem}{\rm
For more than two blocks, the following happens: the space of block diagonal matrices still satisfies condition~\eqref{totgeo}, hence there exists an orthogonal projection on the space of positive definite block diagonal matrices with an arbitrary (fixed) number of blocks, and the corresponding decomposition similar to \eqref{DAD}. However  the space of matrices with zero diagonal blocks do not satisfy condition~\eqref{totgeo} when there are more than two blocks, hence one can not deduce a decomposition of the type \eqref{ADA} in that case.}
\end{rem}

\subsubsection{Decomposition of $SL(2, \mathbb{R})$}
 The group $SL(2, \mathbb{R})$ is the group of isometries of the hyperbolic upper half plane $\mathbb{H}^2$:
$$
\mathbb{H}^2 := \{ z = x + hi, x \in \mathbb{R}, h \in \mathbb{R}^{+*}\} \subset \mathbb{C},
$$
where the action is by homographies.
Mostow's decomposition theorem applied to $$E=\left\{\left(\begin{array}{cc}\alpha & 0\\ 0 & -\alpha\end{array}\right), \alpha \in \mathbb{R}\right\}$$ leads to a decomposition of isometries
$$
SL(2, \mathbb{R}) = SO(2) \times \exp F \times \exp E
$$
where 
\begin{itemize}
\item $SO(2)$ is the stabilizer of $i\in\mathbb{C}$;
\item $\calE = \exp E $ is the space of dilatations $q \mapsto \alpha^2 q$, $\alpha \in \mathbb{R};$
\item $\calF = \exp F = \left\{ \left(\begin{array}{cc} \cosh \beta & \sinh \beta \\ \sinh \beta & \cosh \beta\end{array}\right), \beta \in \mathbb{R}\right\}$ is the space of hyperbolic transformations  with fixed points $-1$ and $1$  in $\mathbb{R}\in\partial\mathbb{H}^2$ and preserving the unit half-circle containing $i$.
\end{itemize}

\subsection{Organization of the paper}

This paper is organized as follows. First we recall the
geometry of the space $\operatorname{SPD}(n)$ of symmetric positive-definite $n\times n$ matrices. We show in particular that it is a
symmetric space of non-positive sectional curvature, homogeneous under the group $\operatorname{GL}(n, \mathbb{R})$ and that the
exponential map defined by the usual power series is  a
diffeomorphism from the space $\operatorname{Sym}(n)$ of
symmetric $n\times n$ real matrices onto $\operatorname{SPD}(n)$.
Moreover we show that the exponential map is the Riemannian
exponential map at the identity with respect to the affine Riemannian metric
on $\operatorname{SPD}(n)$. This is implied by a general result on the
geodesics in locally symmetric spaces that we recall (Proposition \ref{genial}).
This study implies the usual Al-Kashi inequality on the sides of a
geodesic triangle in the non-positively curved space
$\operatorname{SPD}(n)$, and the convexity property of the distance between
two geodesics. In the second subsection, the characterization of the
geodesic subspaces of $\operatorname{SPD}(n)$  given by equation \eqref{totgeo} is proved,  mainly
following arguments from \cite{Mostow}. In subsection \ref{sub3}, the key-step for
the proof of Mostow's decomposition  \eqref{Mostow_dec} is given by the
construction of a non-linear projection from $\operatorname{SPD}(n)$ onto every closed
geodesic subspace. The arguments given here for the existence of
such projection are simpler and more direct then the ones given in
the original paper \cite{Mostow}, and apply to arbitrary
dimension. In the last subsection, we use this projection to prove
the theorem stated above. 

\section{The manifold  $\operatorname{SPD}(n)$ of symmetric positive-definite real  matrices}

Let $\operatorname{M}(n, \mathbb{R})$ be the vector space of real $n\times n$ matrices.
The
group $\operatorname{GL}(n, \mathbb{R})$ is the group of invertible real $n\times n$ matrices, or equivalently the space of real matrices with non-zero determinant:
$$
\begin{array}{l}
\operatorname{GL}(n, \mathbb{R}):= \{g\in \operatorname{M}(n, \mathbb{R}), \det(g) \neq 0\},
\end{array}
$$
The group law given by the multiplication of matrices makes the open set $\operatorname{GL}(n, \mathbb{R})\subset \operatorname{M}(n, \mathbb{R})$ into a Lie-group with Lie-algebra $\operatorname{M}(n, \mathbb{R})$, where the Lie-bracket is given by the commutator of matrices:
\begin{equation}
[A, B] = AB - BA, \forall A,B \in \operatorname{M}(n, \mathbb{R}).
\end{equation}
 The orthogonal group
$\operatorname{O}(n)$ and its Lie algebra
$\mathfrak{so}(n)$ are defined by
$$
\begin{array}{l}
\operatorname{O}(n) := \{ M \in \operatorname{GL}(n, \mathbb{R}), M^T M = M M^T = \operatorname{Id}\}\\ \\
\mathfrak{so}(n) =\{ A \in \operatorname{M}(n, \mathbb{R}), A^T + A = 0\}
\end{array}
$$
The vector space $\operatorname{M}(n, \mathbb{R})$ splits into the direct
sum of $\mathfrak{so}(n)$ and the linear subspace
$\operatorname{Sym}(n)$ of symmetric elements in
$\operatorname{M}(n, \mathbb{R})$
$$
\operatorname{Sym}(n) =\{ M \in \operatorname{M}(n, \mathbb{R}), A^T = A \}.
$$
 The exponential of matrices  defined for all $A$ in $\operatorname{M}(n, \mathbb{R})$ as
\begin{equation}\label{exp}
\exp(A) := \sum_{n = 0}^{+ \infty} \frac{A^{n}}{n!}
\end{equation}
 takes
$\operatorname{Sym}(n)$ to the submanifold $\operatorname{SPD}(n)$ of
$\operatorname{GL}(n, \mathbb{R})$ consisting of symmetric positive-definite
matrices
$$
\exp~:\operatorname{Sym}(n) = \{ A \in \operatorname{M}(n, \mathbb{R}),
A^{T} = A \} \longrightarrow \operatorname{SPD}(n) = \{ A \in
\operatorname{Sym}(n), A^{T}A
> 0 \}.
$$
Note that for $A \in \operatorname{M}(n, \mathbb{R})$, the curve
$\gamma(t):= \exp (t A)$ satisfies 
\begin{equation}\dot{\gamma}(t) =
\left(L_{\gamma(t)}\right)_{*}(A)  = \gamma(t)\cdot A.
\end{equation}
 where $L_{\gamma(t)}$ denotes
left translation by $\gamma(t)$, i.e. multiplication on the left by the matrice $\gamma(t)$, and where the dot denotes the multiplication of matrices. Hence the exponential map defined
by \eqref{exp} is the usual exponential map on the Lie group
$ \operatorname{GL}(n, \mathbb{R})$.
  Let us endowed $\operatorname{SPD}(n) $
with the following Riemannian metric~:
\begin{equation}\label{affine_riemannian}
\langle X, Y\rangle_x := \Tr\left(x^{-1}X x^{-1} Y\right) = \Tr \left(x^{-\frac{1}{2}}X x^{-\frac{1}{2}}\cdot x^{-\frac{1}{2}}Y x^{-\frac{1}{2}}\right) 
\end{equation} 
where $x \in \operatorname{SPD}(n)$,  $X, Y \in T_{x}\operatorname{SPD}(n) $. The identity  matrix $\operatorname{Id}$ belongs to
$\operatorname{SPD}(n)$ and will be our reference point in the following.

\begin{prop}\label{isom}
The left action  of $ \operatorname{GL}(n, \mathbb{R})$ on $\operatorname{SPD}(n) $ defined
by
\begin{equation}\label{theaction}
\begin{array}{cll}
\operatorname{GL}(n, \mathbb{R}) \times \operatorname{SPD}(n)   &\rightarrow&\operatorname{SPD}(n)\\
 (\, g\,,\,x\, ) & \mapsto &  g \cdot x \cdot g^T,
\end{array}
\end{equation}
is a transitive action by isometries (here the dot $\cdot$ denotes the multiplication of matrices).
\end{prop}

\dem \ref{isom}:\\
For every $x$ in $\operatorname{SPD}(n) $, the square root of $x$ is
well-defined and belongs to $\operatorname{SPD}(n)$. In other words there
exists $y$ in $\operatorname{SPD}(n) $ such that $x = y^{2}$. Since $y^T = y$,
one has $x = y\cdot y^T = y\cdot \operatorname{Id} \cdot\, y^T$, and the transitivity follows. To
show that $\operatorname{GL}(n, \mathbb{R})$ acts by isometries, one has to show that the scalar product $\langle \cdot, \cdot \rangle_x$ on the tangent space to the manifold $\operatorname{SPD}(n)$ at $x  = y\cdot \operatorname{Id} \cdot\, y^T$  is related to the scalar product $\langle \cdot, \cdot \rangle_{\textrm{Id}}$ at $\operatorname{Id}$ by
\begin{equation}
\langle X, Y\rangle_{\textrm{Id}} =  \langle y_*X, y_* Y\rangle_{y \cdot \operatorname{Id} \cdot y^T},
\end{equation}
where $y_*: T_{\operatorname{Id}}\operatorname{SPD}(n) \rightarrow T_{{y \cdot \operatorname{Id} \cdot y^T}}\operatorname{SPD}(n)$ denotes the infinitesimal action of $y\in \operatorname{GL}(n, \mathbb{R})$ on tangent vectors to the manifold  $\operatorname{SPD}(n)$. To compute this infinitesimal action, consider a curve $x(t) \in \operatorname{SPD}(n)$ such that $x(0) = x$ and $\dot{x}(t) = X$, and, for a fixed element $g \in \operatorname{GL}(n, \mathbb{R})$, differentiate action \eqref{theaction} to get
$$
g_*X = {\frac{d}{dt}}_{| t = 0}g\cdot x(t)\cdot g^T = g\cdot X \cdot g^T.
$$
It follows that 
\begin{align*}
\langle y_*X, y_* Y\rangle_{y \cdot \operatorname{Id} \cdot y^T} & = \langle y_*X, y_* Y\rangle_{y y^T} = 
\Tr (y y^T)^{-1} (y \cdot X \cdot y^T) (y y^T)^{-1} (y \cdot Y \cdot y^T) \\ &= \Tr (y^T)^{-1} X Y y^T = \Tr (X Y) = \langle X, Y\rangle_{\textrm{Id}}.
\end{align*}
In particular, the scalar product $\langle\cdot, \cdot\rangle_{\textrm{Id}}$ is invariant by the action of the isotropy group of the reference point $\textrm{Id}$, which is exactly the orthogonal group $\operatorname{O}(n)$.
\cqfd

\begin{cor}
The manifold $\operatorname{SPD}(n)$ of symmetric positive-definite real
matrices  is a homogeneous
Riemannian manifold
$$\operatorname{SPD}(n)~=~\operatorname{GL}(n, \mathbb{R})/\operatorname{O}(n).$$
\end{cor}

\begin{thm}\label{diffeo}
The exponential map is a diffeomorphism from
$\operatorname{Sym}(n)$ onto $\operatorname{SPD}(n)$.
\end{thm}

\prooft \ref{diffeo}:\\
Since every $x$ in $\operatorname{SPD}(n)$ admits an orthogonal basis of eigenvectors with positive eigenvalues, the exponential map from $\operatorname{Sym}(n)$ to $\operatorname{SPD}(n)$ is injective and onto, the inverse mapping being given by the logarithm. The differential of the exponential map is recalled in the appendix (Proposition~\ref{log}, equation (\ref{difexp})). By Lemma~\ref{tx}, equation~\eqref{second_difexp}),  it can be written as
$$
d_{X}\exp(Y) = R_{\exp(\frac{X}{2})} L_{\exp(\frac{X}{2})}\tau_{X}(Y),
$$
where $\tau_{X}(Y)$ is the linear isomorphism of $\operatorname{Sym}(n)$ given by equation~\eqref{secondTx}.
Hence $d_{X}\exp$ is continuous and invertible. By the inverse function Theorem, it follows that $\exp$ is a diffeomorphism.  \cqfd

\begin{defe}{\rm
Let $G$ be a  Lie group. An homogeneous space $M = G/K$ is
called \emph{reductive} if the Lie algebra $\mathfrak{g}$ of $G$
splits into a direct sum $\mathfrak{g} = \mathfrak{k} \oplus
\mathfrak{m}$, where $\mathfrak{k}$ is the Lie algebra of $K$, and
$\mathfrak{m}$ an $\textrm{Ad}(K)$-invariant complement. A
reductive homogeneous space is called \emph{locally symmetric} if
the commutation relation $[\mathfrak{m}, \mathfrak{m}] \subset
\mathfrak{k}$ holds.}
\end{defe}

A locally symmetric space is a particular case of a
\emph{naturally reductive} space (see Definition 7.84 page~196 in
\cite{Bes}, Definition 23 page~312 in \cite{ON}, or Proposition
5.2 page~125 in \cite{Ar03} and the definition that follows). In
the finite-dimensional setting, the geodesics of a naturally
reductive space are orbits of one-parameter subgroups of $G$ (see
Proposition 25 page 313 in \cite{ON} for a proof of this fact).
The symmetric case is also treated in Theorem 3.3 page~173 in
\cite{Hel}.  Its infinite-dimensional version has been given in
Example 3.9 in \cite{Ne02c}. For the sake of completeness, we give a proof in the appendix (see Proposition~\ref{genial}) based on the fact that in this case the Levi-Civita connection is the
 homogeneous connection (see \cite{Tum2} for the infinite-dimensional case). We deduce the following theorem.


\begin{thm}\label{formgeo} 
The manifold $\operatorname{SPD}(n)$ is a locally symmetric homogeneous space.
The curve $\gamma(t) := \exp \left(t \log(x)\right)$, $( 0 \leq t
\leq 1)$ is the unique geodesic in
  $\operatorname{SPD}(n)$ joining the identity $~o =\textrm{Id}$ to the element $x
  \in \operatorname{SPD}(n)$. More generally, the geodesic $\gamma_{x,y}(t)$ between any two 
  points $x, y$ in $\operatorname{SPD}(n)$ exists and is unique, and is given by
  \begin{equation}\label{explicite_geo}
  \gamma_{x,y}(t) := 
x^{\frac{1}{2}}\left( \exp t \log \left(x^{-\frac{1}{2}} y
x^{-\frac{1}{2}}\right)\right) x^{\frac{1}{2}}.
\end{equation}
\end{thm}

\prooft \ref{formgeo}:\\
This follows from the same arguments as in \cite{Mostow}, or by
the general result stated in  Proposition \ref{genial}. Indeed the
commutation relation $[\operatorname{Sym}(n), \operatorname{Sym}(n)]
\subset \mathfrak{so}(n)$ implies that $\operatorname{SPD}(n) = \operatorname{GL}(n, \mathbb{R})/\operatorname{O}(n)$ is locally symmetric. It
follows that
$$\gamma(t) := \exp \left(t \log(x)\right) =
\left(\exp \left(\frac{t}{2} \log(x)\right)\right)\cdot \textrm{Id} \cdot \left(\exp \left(\frac{t}{2} \log(x)\right)\right)^T, \qquad
( 0 \leq t \leq 1),
$$
is a geodesic joining $~o = \textrm{Id}$ to $x$. The uniqueness of
this geodesic follows from Proposition~\ref{genial} and Theorem~\ref{diffeo},
since any other geodesic $\gamma_2$ joining $~o = \textrm{Id}$ to $x$
is necessarily of the form 
$$\gamma_{2}(t) = \exp
t\dot{\gamma_{2}}(0) = \left(\exp \left(\frac{t}{2} \dot{\gamma_{2}}(0)\right)\right)\cdot \textrm{Id} \cdot \left(\exp \left(\frac{t}{2} \dot{\gamma_{2}}(0)\right)\right)^T$$
 with velocity $\dot{\gamma_{2}}(0)\in \operatorname{Sym}(n)$. Since $\textrm{GL}(n)$ acts by isometries on $\operatorname{SPD}(n)$, for $x$, $y\in \operatorname{SPD}(n)$, the image of the geodesic $t\mapsto \exp t \log \left(x^{-\frac{1}{2}} y
x^{-\frac{1}{2}}\right)$ by the isometry $u \mapsto x^{\frac{1}{2}} u x^{\frac{1}{2}} = x^{\frac{1}{2}} \cdot u \cdot (x^{\frac{1}{2}})^T$  is itself a geodesic. It follows that the unique geodesic
$\gamma_{x,y}$ joining two points $x$ and $y$ in $\operatorname{SPD}(n)$ is given by
$$\gamma_{x,y}(t) := 
x^{\frac{1}{2}}\left( \exp t \log \left(x^{-\frac{1}{2}} y
x^{-\frac{1}{2}}\right)\right) x^{\frac{1}{2}}.$$ 
\cqfd

\begin{prop}\label{neg}
The manifold $\operatorname{SPD}(n)$ of symmetric positive-definite
real matrices is a Riemannian manifold of non-positive sectional curvature.
\end{prop}

\dem \ref{neg}:\\
The Levi-Civita connection is defined
by Koszul formula (see for instance Theorem 3.1 page~54 in
\cite{Ar03} where this formula is recalled) and coincides with the homogeneous connection since $\operatorname{SPD}(n)$ is locally symmetric and the Riemannian metric is $G$-invariant. The sectional
curvature $K_{\textrm{Id}}$ at the reference point $ \textrm{Id}$ is given by~:
$$
K_{\textrm{Id}}(X, Y) := \frac{\langle R_{X, Y}X, Y\rangle_\textrm{Id}}{\langle X,
X\rangle_\textrm{Id}\langle Y, Y\rangle_\textrm{Id} -\langle X, Y\rangle_\textrm{Id}^{2}} = \frac{\langle \left[ [X, Y],
X\right], Y \rangle_\textrm{Id}}{\langle X, X\rangle_\textrm{Id} \langle Y, Y\rangle_\textrm{Id} -\langle X, Y\rangle_\textrm{Id}^{2}},
$$
for all $X$, $Y$ in $T_{\textrm{Id}}\operatorname{SPD}(n) = \operatorname{Sym}(n)$, where we have used formula~\eqref{curvature_symmetric} for the curvature (see also
Proposition 7.72 page~193 in \cite{Bes}, or Proposition 6.5
page~92 in \cite{Ar03}). Now the sign of the sectional curvature
of the $2$-plane generated by $X$ and $Y$ is the sign of
$\langle \left[ [X, Y], X\right], Y \rangle_{\textrm{Id}}$. One has~:
$$
\begin{array}{ll}
\langle \left[ [X, Y], X\right], Y \rangle_\textrm{Id}  & = \Tr \left[ [X, Y],
X\right] Y = \Tr \left([X, Y] X Y - X [X, Y] Y \right)\\& = \Tr
[X, Y] [X, Y] = - \Tr [X, Y]^{T} [X, Y] \leq 0,
\end{array}
$$
where the last identity following from the fact that  $[X, Y]$
belongs to $[\operatorname{Sym}(n), \operatorname{Sym}(n)] \subset
\mathfrak{so}(n)$.
 \cqfd
\\
\\

The following two Lemmas are standard results in the geometry of
non-positively curved spaces.

\begin{lem}\label{angle}
The Riemannian angle between two paths $f$ and $g$ intersecting at
$o = \textrm{Id}$ is equal to the Euclidian angle between the two paths
$\log(f)$ and $\log(g)$ at $0$. Moreover, in any geodesic triangle
$ABC$ in $\operatorname{SPD}(n)$,
\begin{equation}\label{Al-Kashi}
c^{2} \geq a^{2} + b^{2} - 2 a b \cos \widehat{ACB},
\end{equation}
where $a, b, c$ are the lengths of the sides opposite to $A, B, C$ respectively,
and $ \widehat{ACB}$ the angle at $C$.
\end{lem}

\demlem \ref{angle}:\\ This follows from the same arguments as in
\cite{Mostow}. The Al-Kashi inequality \eqref{Al-Kashi} is also a direct
consequence of  Corollary 13.2 in \cite{Hel} page 73, since Lemma
\ref{formgeo} implies that $\operatorname{SPD}(n)$ is a minimizing convex
normal ball.\cq


\begin{lem}\label{convex}
Let $\gamma_{1}(t)$ and $\gamma_{2}(t)$ be two constant speed
geodesics in $\operatorname{SPD}(n)$. Then the distance in $\operatorname{SPD}(n)$
between $\gamma_{1}(t)$ and $\gamma_{2}(t)$ is a convex function
of $~t$.
\end{lem} \cq

\demlem \ref{convex}:\\
This follows from the fact that a non-positively curved Riemannian manifold is a CAT(0) space.
\\
\\

Let us recall the following definition.
\begin{defe}{\rm
A Riemannian manifold $\mathcal{P}$ is called
\emph{symmetric} if for every $p$ in $\mathcal{P}$, there exists a
globally defined isometry $s_{p}$ which fixes $p$ and such that
the differential of $s_{p}$ at $p$ is $-\id$. The transformation $s_{p}$ is called a \textit{global symmetry with respect to} $p$.}
\end{defe}

\begin{prop}\label{symneg}
The manifold $\operatorname{SPD}(n)$  is a globally symmetric homogeneous Riemannian
manifold. The globally defined symmetry with respect to $x\in\operatorname{SPD}(n)$ is
\begin{equation}\label{symmetry}
s_{x}(y) = x^{\frac{1}{2}}(x^{-\frac{1}{2}} y x^{-\frac{1}{2}})^{-1} x^{\frac{1}{2}} = x~y^{-1}\,x.
\end{equation}
\end{prop}

\dem \ref{symneg}:\\
Consider the inversion $\operatorname{Inv}\,: \operatorname{SPD}(n) \rightarrow \operatorname{SPD}(n)$
defined by $\operatorname{Inv}(x) = x^{-1}$. An easy computation shows that the
differential $T_x\operatorname{Inv}$ of $\operatorname{Inv}$ at $x \in \operatorname{SPD}(n)$ takes a tangent
vector $X \in T_{x}\operatorname{SPD}(n)$ to $T_x\operatorname{Inv}(X) = -x^{-1} X x^{-1}$.
One therefore has~:
$$
\begin{array}{ll}
\langle  T_x\operatorname{Inv}(X), T_x\operatorname{Inv}(Y) \rangle_{\operatorname{Inv}(x)} & = \langle
-x^{-1} X x^{-1}, -x^{-1} Y x^{-1}\rangle_{x^{-1}} = \Tr x(-x^{-1} X
x^{-1}) x (-x^{-1} Y x^{-1})\\
& = \Tr X x^{-1} Y x^{-1} = \Tr x^{-1} X x^{-1} Y = \langle X, Y\rangle_{x},
\end{array}
$$
hence $\operatorname{Inv}$ is an isometry of $\operatorname{SPD}(n)$. Since $\operatorname{Inv}$ fixes the reference point 
 $\textrm{Id}\in \operatorname{SPD}(n)$ and its differential at the reference point is
$-\textrm{Id}$ on $T_{\textrm{Id}}\operatorname{SPD}(n) = \operatorname{Sym}(n)$, it follows
that $\operatorname{Inv}$ is a global symmetry with respect to $\textrm{Id}$. From
the transitive action of $\operatorname{GL}(n, \mathbb{R})$ one can define a globally defined symmetry with respect to any $x\in\operatorname{SPD}(n)$. Indeed note that $s_x$ is the composition of the isometry $\iota_x: y \mapsto x^{-\frac{1}{2}}\cdot y \cdot x^{-\frac{1}{2}}$, with the inversion $\operatorname{Inv}$ followed by the isometry $(\iota_x)^{-1}: y \mapsto x^{\frac{1}{2}}\cdot y \cdot x^{\frac{1}{2}}$. Hence $s_x$ is an isometry, with fixed point $x$.  Its differential at $x$ is given by $$T(\iota_x)^{-1} \circ T_{\textrm{Id}}\operatorname{Inv}\circ T_x\iota_x = T_{\textrm{Id}}(\iota_x)^{-1} \circ (-\textrm{Id}) \circ T_x\iota_x = -\textrm{Id}.$$
Whence
$\operatorname{SPD}(n)$ is a (globally) symmetric homogeneous Riemannian manifold.
\cqfd

\section{Totally geodesic submanifolds of the manifold $\operatorname{SPD}(n)$}

\begin{defe}
{\rm
A submanifold $\mathcal{N}$ of a Riemannian manifold $\mathcal{M}$ is called \textit{totally geodesic}   if any geodesic on $\mathcal{N}$ for the induced Riemannian metric is also a geodesic on $\mathcal{M}$. In particular, the property of being totally geodesic is a local property of the submanifold $\mathcal{N}$ inside the ambient manifold $\mathcal{M}$.}
\end{defe}
\begin{defe}
{\rm
A submanifold $\mathcal{N}$ of a Riemannian manifold $\mathcal{M}$ is called \textit{geodesically convex}   if there exists a unique minimizing geodesic of $\mathcal{M}$ connecting two points in $\mathcal{N}$ and entirely contained in $\mathcal{N}$.  In particular, the property of being geodesically convex depends on global behavior of geodesics.

}
\end{defe}
The following Theorem is based on \cite{Mostow}. We give below a self-contained proof based on the following Lemma.

\begin{thm}[\cite{Mostow}]\label{equivgeospace}
Let $\operatorname{E}$ be a linear subspace of $\operatorname{Sym}(n)$. The
following assertions are equivalent~:
\begin{enumerate}
\item \label{m1} $[ X, [X, Y]] \in E$ for all  $X, Y \in \operatorname{E}$,
\item\label{m2} $e\cdot  f \cdot e \in \calE := \exp(\operatorname{E})$ for all $e, f \in \calE = \exp(\operatorname{E})$,
\item \label{m5} $\calE := \exp(\operatorname{E})$ is a geodesically convex submanifold of
  $\operatorname{SPD}(n)$.
\item \label{m3} $\calE := \exp(\operatorname{E})$ is a totally geodesic submanifold of
  $\operatorname{SPD}(n)$.
  \item \label{m4} $\operatorname{E}$ is a Lie triple system, i.e. for all $X, Y, Z\in \operatorname{E}$, $[X, [Y, Z]]\in E$.
\end{enumerate}
\end{thm}

\begin{lem}\label{techu}:\\
Let $f\in\calE = \exp(\operatorname{E})$,
$Y\in \operatorname{E}$, and consider the differentiable curve  in
$\operatorname{Sym}(n)$ defined by
$$
X(t) := \log(\exp tY\cdot f \cdot \exp tY).
$$
Then 
\begin{equation}
\dot{X}(t) = \ad(X(t))\coth (\ad(X(t)/2))(Y).
\end{equation}
\end{lem}

\demlem \ref{techu}:\\
One has 
$$
{\frac{d}{dt}}_{|t = t_{0}} \exp tY = Y\cdot \exp t_0Y = \exp t_0Y \cdot Y,
$$
hence the differential of the right hand side of equation~\eqref{lr} is
$$
{\frac{d}{dt}}_{|t = t_{0}} \exp tY\cdot  f \cdot \exp tY = Y \exp X(t_0) + \exp X(t_0) Y.
$$
On the other hand, the differential of the left hand side of equation~\eqref{lr} is
$$
{\frac{d}{dt}}_{|t = t_{0}} \exp X(t) =  (d_{X(t_0)}\exp) (\dot{X}(t_0)).
$$
Hence $X(t)$ satisfies the following differential equation~:
$$
\dot{X}(t) = (d_{X(t)}\exp)^{-1} (Y\cdot \exp X(t) + \exp X(t)\cdot Y).
$$
By Lemma~\ref{tx}, equation~\eqref{second_difexp2},
$$
(d_{X}\exp)^{-1}(Z) = \frac{\ad(X/2)}{\sinh (\ad(X/2))}\left( \exp\left(-\frac{X}{2}\right)\cdot Z \cdot \exp\left(-\frac{X}{2}\right)\right)
$$
It follows that 
\begin{align*}
\dot{X}(t) &= \frac{\ad(X(t)/2)}{\sinh (\ad(X(t)/2))}\left( \exp\left(-\frac{X(t)}{2}\right)\cdot \left(Y\cdot \exp X(t) + \exp X(t)\cdot Y\right) \cdot \exp\left(-\frac{X(t)}{2}\right)\right)\\
&= \frac{\ad(X(t)/2)}{\sinh (\ad(X(t)/2))}\left( \exp\left(-\frac{X(t)}{2}\right)\cdot Y\cdot \exp \left(\frac{X(t)}{2}\right) + \exp \left(\frac{X(t)}{2}\right)\cdot Y \cdot \exp\left(-\frac{X(t)}{2}\right)\right)\\
&= \frac{\ad(X(t)/2)}{\sinh (\ad(X(t)/2))}\left( \textrm{Ad}\left(\exp\left(-\frac{X(t)}{2}\right)\right)(Y) + \textrm{Ad}\left(\exp \left(\frac{X(t)}{2}\right)\right)(Y)\right)\\
&= \frac{\ad(X(t)/2)}{\sinh (\ad(X(t)/2))}\left( \exp\left(\ad\left(-\frac{X(t)}{2}\right)\right)(Y) + \exp\left(\ad\left(\frac{X(t)}{2}\right)\right)(Y)\right)\\
&= \frac{\ad(X(t))}{\sinh (\ad(X(t)/2))}\left(\cosh(\ad(X(t)/2))\right)(Y)\\
&= \ad(X(t))\coth (\ad(X(t)/2))(Y)\\
\end{align*}
\cq

\demthm \ref{equivgeospace}:

 $\ref{m1} \Rightarrow \ref{m2}:$ Let $f\in\calE$ and
$Y\in \operatorname{E}$. We will  show that 
\begin{equation}\label{lr}
\exp X(t) = \exp tY\cdot  f \cdot \exp tY 
\end{equation}  
belongs to $\calE = \exp(\operatorname{E})$ for all $t \in \mathbb{R}$. By Lemma~\ref{techu}, 
\begin{equation*}
\dot{X}(t) = \ad(X(t))\coth (\ad(X(t)/2))(Y).
\end{equation*}
Consider the vector field $W$ on $\operatorname{Sym}(n)$ defined by
$$
W(X) := \ad(X)\coth (\ad(X/2))(Y), X\in \operatorname{Sym}(n).
$$
Note that  only even powers of $\ad(X)$ are involved in the operator
$\ad(X)\coth (\ad(X/2))$. Since $\ad(X)\circ \ad(X)(Y)\in \operatorname{E}$ for $X, Y \in E$, the vector field $W(X)$ is tangent to $\operatorname{E}$ for every $X\in \operatorname{E}$. It follows that an integral curve $X(t)$ of $W$ starting at $X(0)\in \operatorname{E}$ stays in $\operatorname{E}$. Moreover the flow of this vector field
 is defined for all  $t \in
\R$. Thus setting $t = 1$ and $Y = \log e$ with $e \in \calE = \exp \operatorname{E}$,
give the result $e\cdot f \cdot e \in \calE$.

 $\ref{m2} \Rightarrow \ref{m5}:$ Suppose that for all $e$ and
$f$ in $\calE = \exp \operatorname{E}$, the product  $e\cdot f \cdot e$ belongs to $\calE = \exp \operatorname{E}$. Consider $f\in \calE = \exp \operatorname{E}$. By Lemma \ref{formgeo}, the geodesic of $\operatorname{SPD}(n)$ joining $\textrm{Id}\in\calE = \exp \operatorname{E}$ to $f$ is $t\mapsto \exp(t \log(f))$ hence lies entirely in the manifold $\calE = \exp \operatorname{E}$.  In particular, its mid-point $\exp(\frac{1}{2}\log(f)) =: f^{\frac{1}{2}}$ belongs to $\calE = \exp \operatorname{E}$. In other words, the square roots of elements in $\calE = \exp \operatorname{E}$ belong to $\calE = \exp \operatorname{E}$. Consider now $x$ and $y$ two elements in $\calE = \exp \operatorname{E}$. By hypothesis, $x^{-\frac{1}{2}}\cdot y \cdot x^{-\frac{1}{2}}$ belongs to $\calE = \exp \operatorname{E}$, hence $\log x^{-\frac{1}{2}}\cdot y \cdot x^{-\frac{1}{2}}$ belongs to $\operatorname{E}$ and the geodesic 
$$
t \mapsto \exp(t \log x^{-\frac{1}{2}}\cdot y \cdot x^{-\frac{1}{2}})
$$
joining $\textrm{Id}$ to  $x^{-\frac{1}{2}}\cdot y \cdot x^{-\frac{1}{2}}$ lies entirely in $\calE = \exp \operatorname{E}$. 
 Since $x\in \calE = \exp \operatorname{E}$, by hypothesis the isometry $z \mapsto x^{\frac{1}{2}}\cdot z \cdot x^{\frac{1}{2}}$ preserves $\calE = \exp \operatorname{E}$ and sends the previous geodesic to the unique geodesic of $\operatorname{SPD}(n)$ joining $x$ to $y$, which therefore lies completely in $\calE = \exp \operatorname{E}$. Consequently, the space  $\calE = \exp \operatorname{E}$
is geodesically convex in $\operatorname{SPD}(n)$.

$\ref{m5} \Rightarrow \ref{m3}:$ Suppose that $\calE = \exp \operatorname{E}$
is geodesically convex in $\operatorname{SPD}(n)$. By the uniqueness of the geodesic joining two points  in $\operatorname{SPD}(n)$ proved in Theorem~\ref{formgeo}, the geodesic of $\operatorname{SPD}(n)$ connecting two points $x$ and $y$ belonging to $\calE = \exp \operatorname{E}$ minimises the distance between $x$ and $y$ in $\operatorname{SPD}(n)$ and by hypothesis lies entirely in  $\calE = \exp \operatorname{E}$. Therefore it is also a geodesic of $\calE = \exp \operatorname{E}$. Now consider an arbitrary geodesic $t \mapsto \gamma(t)$  in  $\calE = \exp \operatorname{E}$, with $t$ in some open interval $I$, and let us show that it is a geodesic in $\operatorname{SPD}(n)$. For every $t_0\in I$, there exists an $\varepsilon>0$ such that $(t_0-\varepsilon, t_0+\varepsilon)\ni t\mapsto \gamma(t)$ is the unique geodesic in $\calE = \exp \operatorname{E}$ between $\gamma(t_0-\varepsilon)$ and $\gamma(t_0+\varepsilon)$. By uniqueness,  it is the geodesic segment in $\operatorname{SPD}(n)$ connecting $\gamma(t_0-\varepsilon)$ and $\gamma(t_0+\varepsilon)$. Therefore $t \mapsto \gamma(t)$ is a geodesic in $\operatorname{SPD}(n)$.

 $\ref{m3} \Rightarrow \ref{m2}:$ Suppose that  $\calE = \exp \operatorname{E}$ is a
totally geodesic submanifold of $\operatorname{SPD}(n)$. Let us consider
the symmetry $s_{x}$ with respect to $x \in \operatorname{SPD}(n)$ defined
from $\operatorname{SPD}(n)$ to $\operatorname{SPD}(n)$ by $s_{x}: y \mapsto x y^{-1}
x$ (see Proposition~\ref{symneg}).  Every geodesic of the form $t
\mapsto x^{\frac{1}{2}}\exp(tY) x^{\frac{1}{2}}$ is mapped to $t
\mapsto x^{\frac{1}{2}}\exp(-tY) x^{\frac{1}{2}}$ by $s_{x}$. It
follows that every geodesic containing   $x$ is stable under
$s_{x}$. Consequently  $s_{x}(\exp \operatorname{E}) \subset \exp \operatorname{E}$ for every $x
\in \exp \operatorname{E}$. In particular, for $e$ and $f$ in $\calE = \exp \operatorname{E}$, 
$$
e\cdot f \cdot e = e^2\cdot \left(e\cdot f^{-1} \cdot e\right)^{-1}\cdot e^2 = s_{e^2}\left(s_e(f)\right)
$$
belongs to $\calE = \exp \operatorname{E}$.

 $\ref{m2} \Rightarrow \ref{m1}$: Suppose that $e\cdot  f \cdot e \in \calE := \exp(\operatorname{E})$ for all $e, f \in \calE = \exp(\operatorname{E})$. Given $X$ and $Y$ in $\operatorname{E}$, consider the smooth function in
$\operatorname{Sym}(n)$ defined by
$$
Z(s, t) = \log(\exp tY\cdot \exp sX \cdot\exp tY).
$$
In particular,
$$
Z(s, 0) = \log(\exp sX) = s X.
$$
By hypothesis,  $Z(s, t)$ belongs to $\operatorname{E}$ for all $(s, t) \in \R^2$.  Since $\operatorname{E}$ is a vector space, all the derivatives of the function $(s, t) \mapsto Z(s, t)$ belong to $\operatorname{E}$ for all $(s, t) \in \R^2$. By Lemma~\ref{techu},  $$\frac{\partial}{\partial t}Z(s, t) =
\ad(Z(s, t))\coth\left(\ad(Z(s, t)/2)\right)(Y).$$ 
In particular
$${\frac{\partial}{\partial t}}_{| t = 0}Z(s, t) =
\ad(Z(s, 0))\coth\left(\ad(Z(s, 0)/2)\right)(Y) = \ad(sX)\coth\left(\ad(sX/2)\right)(Y)$$
belongs to $\operatorname{E}$ for all $s\in \R$.
Consider the Taylor series of $\frac{u}{2} \coth\left(\frac{u}{2}\right)$ in the neighborhood of $u = 0$:
$$
\frac{u}{2} \coth\left(\frac{u}{2}\right) = 1 +  \frac{u^2}{12} + u^2\varepsilon(u), \textrm{ with }\lim_{u \mapsto 0} \varepsilon(u) = 0.
$$
It follows that
$$
{\frac{\partial}{\partial t}}_{| t = 0}Z(s, t) = 2\left(Y + \frac{1}{12}\ad(sX)\circ \ad(sX)(Y) + s^2 W\right)\in \operatorname{E} \textrm{ with }\lim_{s \mapsto 0} W = 0.
$$
Therefore
$$
\lim_{s \rightarrow 0}\frac{{\frac{\partial}{\partial t}}_{| t = 0}Z(s, t) - 2 Y }{s^2} = \frac{1}{6} [X, [X, Y]]\in \operatorname{E}.
$$

$\ref{m1} \Rightarrow \ref{m4}$: Consider  $X, Y, Z$ in $\operatorname{E}$. Developping $[X + Y, [X + Y, Z]]\in \operatorname{E}$ we get 
$$
[X, [Y, Z]] + [Y, [X, Z]]\in \operatorname{E}.
$$ 
By Jacobi identity,
$$
[ Y, [X, Z]] = [[Y, X], Z] + [X, [Y, Z]] \in \operatorname{E},
$$
hence 
$$
A := [[Y, X], Z] + 2 [X, [Y, Z]] \in \operatorname{E}.
$$
Interchanging $X$ and $Z$ we also have
$$
B := [[Y, Z], X] + 2 [Z, [Y, X]] \in \operatorname{E}.
$$
It follows that 
$$
2 A + B = 2 [[Y, X], Z] + 4 [X, [Y, Z]]  + [[Y, Z], X] + 2 [Z, [Y, X]] = 3 [X, [Y, Z]] \in \operatorname{E}.
$$

$\ref{m4} \Rightarrow \ref{m1}$ Obvious by taking $X = Y$.
\cqfdt

\section{Orthogonal projection on a totally geodesic
submanifold}\label{sub3}

The proof of Mostow's decomposition theorem given in \cite{Mostow}
is based on the existence of an orthogonal projection from
$\operatorname{SPD}(n)$ onto a totally geodesic submanifold $\calE := \exp \operatorname{E}$ which follows from compactness
arguments. 
Here we use the completeness of $\calE := \exp \operatorname{E}$ to obtain this projection. Recall from \cite{Bhatia} the following Proposition (we sketch the proof for the convenience of the reader).

\begin{prop}\label{complete}
The manifold $\operatorname{SPD}(n)$ endowed with the distance induced from the Riemannian metric \eqref{affine_riemannian} is a complete metric space. For any linear subspace $\operatorname{E}$ of $\operatorname{Sym}(n)$, the manifold $\calE := \exp \operatorname{E}$, endowed with the induced Riemannian metric, is closed in $\operatorname{SPD}(n)$ hence also a complete metric space.
\end{prop}

\dem \ref{complete}
By \eqref{secondTx}, the exponential map $\exp: \operatorname{Sym}(n) \rightarrow \operatorname{SPD}(n)$ increases distances, hence a Cauchy sequence in $\operatorname{SPD}(n)$ is the image by the exponential map of a Cauchy sequence in $ \operatorname{Sym}(n)$ which is complete. By Theorem~\ref{diffeo}, it follows that $\operatorname{SPD}(n)$  is complete. The remainder follows from the fact that a linear subspace $\operatorname E$ is closed in $\operatorname{Sym}(n)$.
\cqfd

\begin{thm}\label{projection}
Let $\operatorname{E}$ be a linear subspace of $\operatorname{Sym}(n)$ such that $[ X, [X, Y]] \in \operatorname{E}$, for all $X, Y \in \operatorname{E}$.
Then there exists a continuous orthogonal projection  from $\operatorname{SPD}(n)$
onto the totally geodesic submanifold $\calE := \exp \operatorname{E}$, i.e. a continuous map $\pi: \operatorname{SPD}(n) \rightarrow \calE$  satisfying
$$\textrm{dist}(x, \exp \operatorname{E}) = \textrm{dist}(x, \pi(x))$$ and such
that the geodesic joining $x$ to $\pi(x)$ is orthogonal to every
geodesic starting at $\pi(x)$ and included in $\calE := \exp \operatorname{E}$.
\end{thm}

\begin{figure}
\begin{center}
\includegraphics[width = 0.8 \textwidth]{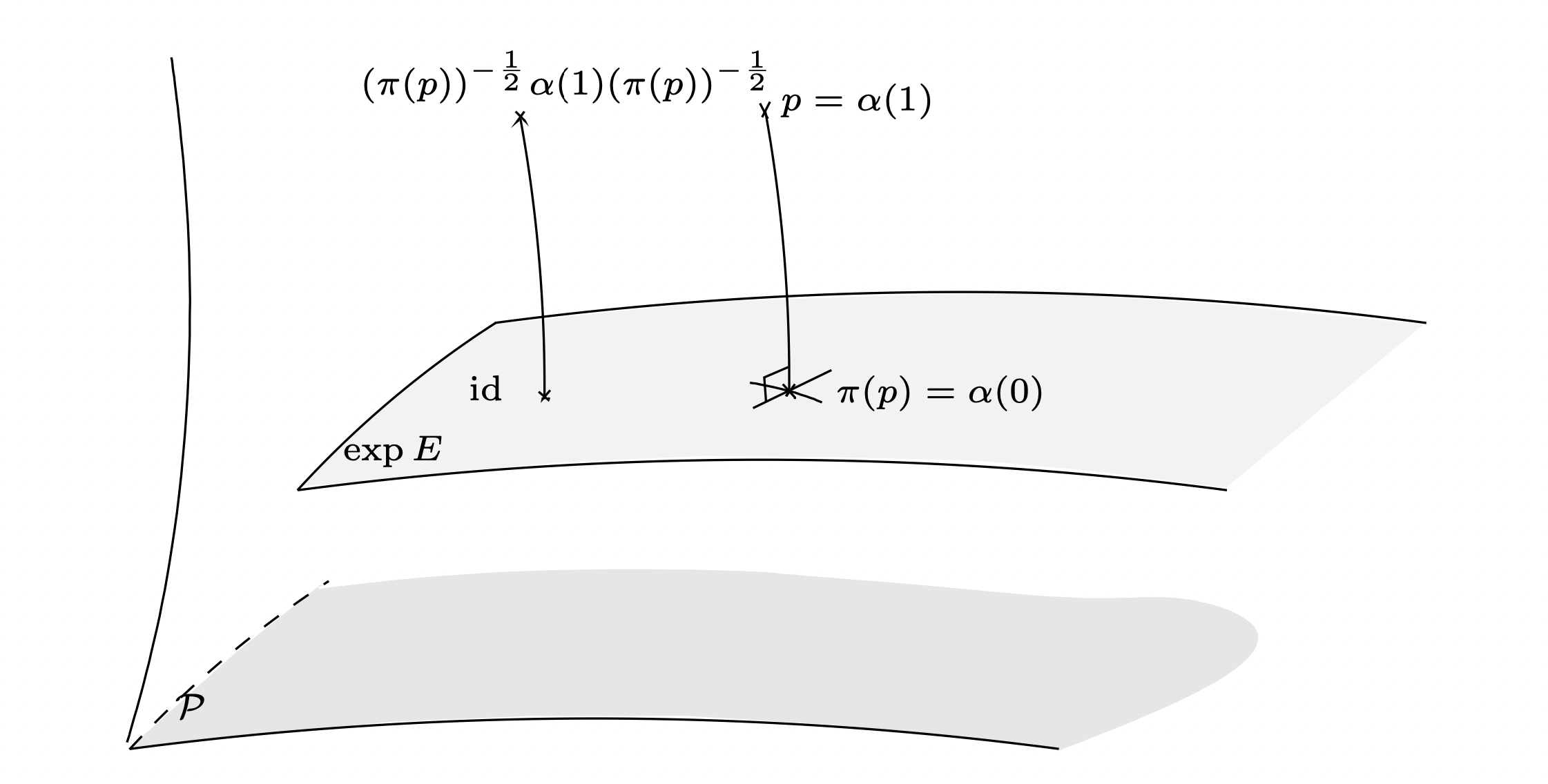}
\caption{Illustration of the orthogonal projection $\pi(p)$ of a point $p\in SPD(n)$ onto a totally geodesic submanifold $\exp E$ and the translation of the fiber over $\pi(p)$ by the action of $\pi(p)^{-\frac{1}{2}}.$}
\end{center}
\end{figure}

\begin{figure}
\begin{center}
\includegraphics[width = 0.8 \textwidth]{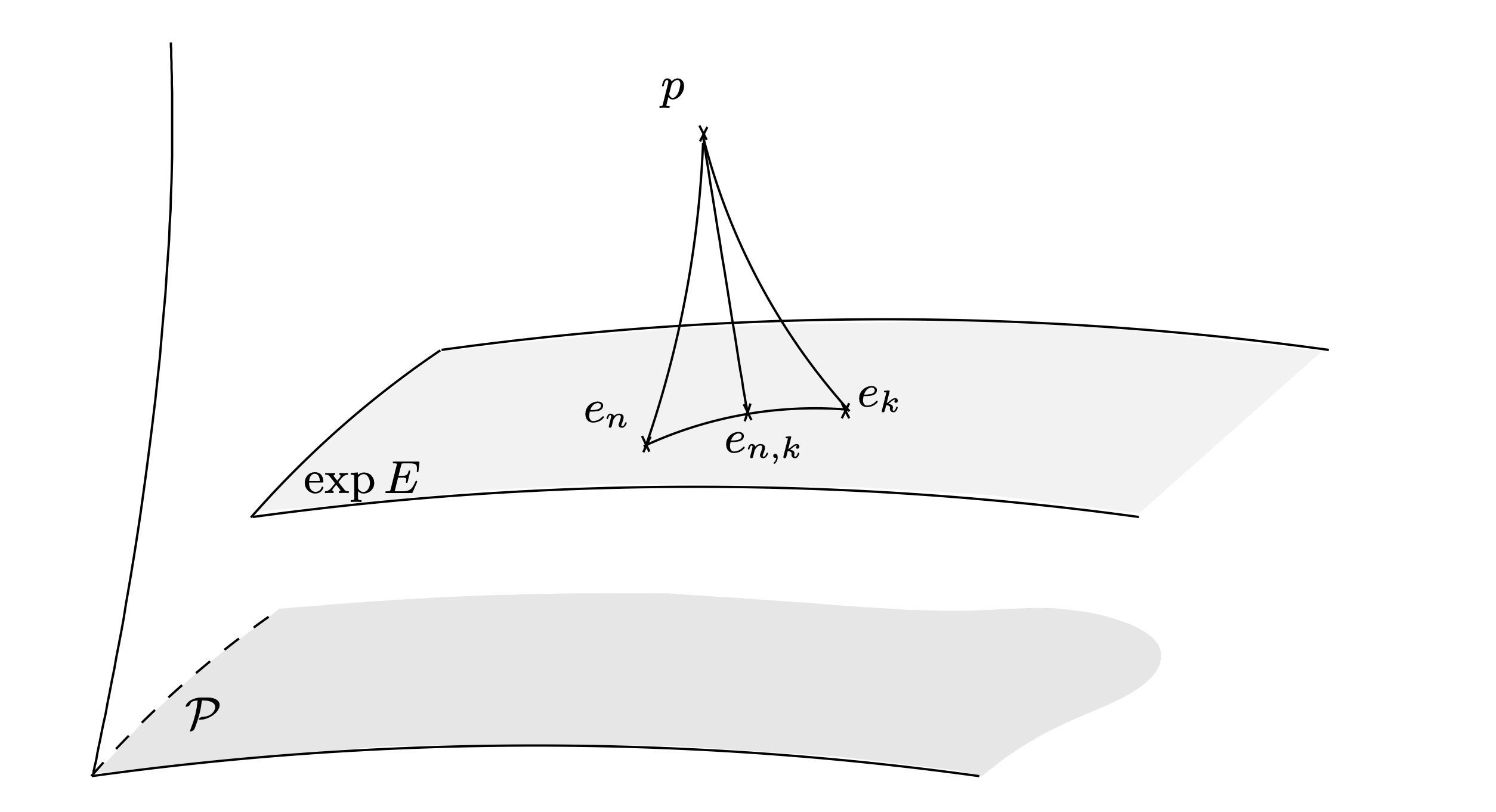}
\caption{Illustration of geodesic triangles used in  the proof of Therorem~\ref{projection}. }
\end{center}
\end{figure}

\prooft \ref{projection}:\\
 Let $x$ be a element of $\operatorname{SPD}(n)$. Denote by $\delta$ the
distance between $x$ and $\calE := \exp \operatorname{E}$ in $\operatorname{SPD}(n)$ and let
$\{e_{n} \}_{n \in \N}$ be a minimizing sequence in $\calE := \exp \operatorname{E}$ thus that
$$
\textrm{dist}(x, e_{n})^{2} \leq \delta^{2} + \frac{1}{n}.
$$
Let us show that $\{ e_{n} \}_{n \in \N}$ is a Cauchy sequence in
$\calE := \exp \operatorname{E}$. For this purpose, consider for $k > n$ the geodesic
$\gamma(t)$ joining $e_{n} =: \gamma(0)\in \calE$ to $e_{k}:= \gamma(1)\in \calE$.
This geodesic  lies in $\calE$ since $\calE$ is a totally geodesic
submanifold of $\operatorname{SPD}(n)$, and is of the form:
$$
\gamma(t) = e_{n}^{\frac{1}{2}} \exp(t H) e_{n}^{\frac{1}{2}},
$$
where $H$ belongs  to $\operatorname{E}$. Denote by $e_{n, k}$ the middle of the
geodesic joining $e_{n}$ to $e_{k}$, i.e. $$e_{n, k} =
e_{n}^{\frac{1}{2}} \exp \left(\frac{H}{2} \right) e_{n}^{\frac{1}{2}}.$$ 
By
lemma \ref{angle} applied to the geodesic triangle joining $x$,
$e_{n}$ and $e_{n, k}$, we have:
\begin{equation}\label{eq1}
\textrm{dist}(x, e_{n})^{2} \geq \textrm{dist}(e_{n}, e_{n,
k})^{2} + \textrm{dist}(e_{n, k}, x)^{2} - 2 \textrm{dist}(e_{n},
e_{n, k}) \textrm{dist}(e_{n, k}, x) \cos \widehat{e_{n} e_{n, k}
x}.
\end{equation}
On the other hand, lemma \ref{angle} applied to the geodesic
triangle joining $x$, $e_{k}$ and $e_{n, k}$ gives:
\begin{equation}\label{eq2}
\textrm{dist}(x, e_{k})^{2} \geq \textrm{dist}(e_{k}, e_{n,
k})^{2} + \textrm{dist}(e_{n, k}, x)^{2} - 2 \textrm{dist}(e_{k},
e_{n, k}) \textrm{dist}(e_{n, k}, x) \cos \widehat{e_{k} e_{n, k}
x}.
\end{equation}
By definition of $e_{n, k}$ we have: $\textrm{dist}(e_{k}, e_{n,
k}) = \textrm{dist}(e_{n}, e_{n,
  k})$. Moreover since the geodesic $\gamma$ is a smooth curve, the sum of the following angles is the flat one:
$$
\widehat{e_{k} e_{n, k} x} + \widehat{e_{n} e_{n, k} x} =\pi,
$$
and $\cos \widehat{e_{k} e_{n, k} x} = - \cos \widehat{e_{n} e_{n,
k}
  x}.$
Summing  inequalities \eqref{eq1} and \eqref{eq2}, we obtain:
$$
\textrm{dist}(x, e_{n})^{2} + \textrm{dist}(x, e_{k})^{2} \geq 2
\textrm{dist}(e_{k}, e_{n, k})^{2} + 2 \textrm{dist}(e_{n, k},
x)^{2}.
$$
It follows that:
\begin{align*}
\textrm{dist}(e_{k}, e_{n, k})^{2} & \leq
\frac{1}{2}\left(\textrm{dist}(x, e_{n})^{2} + \textrm{dist}(x,
e_{k})^{2}\right) - \textrm{dist}(e_{n, k},
x)^{2}\\
& \leq \frac{1}{2}\left(\delta^{2} + \frac{1}{n} + \delta^{2} +
\frac{1}{k}\right) -
\delta^{2}\\
& \leq \frac{1}{2}\left(\frac{1}{n} + \frac{1}{k}\right).
\end{align*}
This yields that $$\textrm{dist}(e_{n}, e_{k}) \leq
\sqrt{2}\left(\frac{1}{n} + \frac{1}{k}\right)^{\frac{1}{2}}.$$ Consequently
$\{e_{n}\}_{n
  \in \N}$ is a Cauchy sequence in $\calE = \exp \operatorname{E}$.
By Proposition~\ref{complete}, $\calE$ is a complete metric space hence the sequence $\{ e_{n} \}_{n \in \N} $
converges to a element $\pi(x)$ in $\calE$ satisfying:
$$
 \textrm{dist}(x, \pi(x)) = \textrm{dist}(x, \calE).
$$

 Denote by $\alpha(t)$ the constant speed geodesic
 which satisfies $\alpha(0) = \pi(x)$ and $\alpha(1) = x$. By
uniqueness of the geodesic joining two points in $\operatorname{SPD}(n)$, it follows that the
length of $\alpha $ is $\textrm{dist}(x, \calE)$. The map $$y
\mapsto (\pi(x))^{-\frac{1}{2}} y (\pi(x))^{-\frac{1}{2}}$$ being
an isometry of $\operatorname{SPD}(n)$ which preserves $\calE$, the curve $t \mapsto (\pi(x))^{-\frac{1}{2}} \alpha(t)
(\pi(x))^{-\frac{1}{2}}$ is a geodesic whose length is the
distance between $$(\pi(x))^{-\frac{1}{2}} \alpha(1)
(\pi(x))^{-\frac{1}{2}} = (\pi(x))^{-\frac{1}{2}} x
(\pi(x))^{-\frac{1}{2}}$$ and $\calE$. Therefore the projection of
$(\pi(x))^{-\frac{1}{2}} x (\pi(x))^{-\frac{1}{2}}$ onto $\calE$
is $(\pi(x))^{-\frac{1}{2}} \alpha(0) (\pi(x))^{-\frac{1}{2}} = \textrm{Id}$. From lemma \ref{formgeo} it follows that:
$$
(\pi(x))^{-\frac{1}{2}} \alpha(t) (\pi(x))^{-\frac{1}{2}} = \exp t
V,
$$
for some $V$ in $\operatorname{Sym}(n)$. Since the length of $t
\mapsto \exp tV$ is $\|V\|$, $V$  is in $F$ and
$(\pi(x))^{-\frac{1}{2}} x (\pi(x))^{-\frac{1}{2}}$ is in $\exp
F$. Since $E \perp F$, by lemma \ref{angle},
$(\pi(x))^{-\frac{1}{2}} \alpha (\pi(x))^{-\frac{1}{2}}$ is
orthogonal at the identity to every curve starting at the identity
and contained in $\exp E$. Therefore $\alpha$ is orthogonal at
$\pi(x)$ to every curve starting at $\pi(x)$ and contained in
$\exp E$.

 To show that $\pi$ is continuous, denote by $\gamma(t)$
 (resp. $\alpha(t)$) the geodesic
joining a  point $x_{1}$ (resp. $x_{2}$) in $\operatorname{SPD}(n)$ to its
projection onto $\exp E$, with $\gamma(0) = \pi(x_{1})$ (resp.
$\alpha(0) = \pi(x_{2})$ ) and $\gamma(1) = x_{1}$ (resp.
$\alpha(1) = x_{2}$). By the negative curvature property stated in
Lemma \ref{convex}, the map $t \mapsto \textrm{dist}(\gamma(t),
\alpha(t))$ is convex. Since, for $t = 0$, $\gamma(t)$ and
$\alpha(t)$ are orthogonal to the geodesic joining $\pi(x_{1})$
and $\pi(x_{2})$ (which is contained in $\calE$), the minimum of the distance between $\gamma(t)$
and $\alpha(t)$ is reached for $t = 0$, and $\textrm{dist}(x_{1},
x_{2}) \geq \textrm{dist}(\pi(x_{1}), \pi(x_{2}))$. \cqfdt

\section{Mostow's decomposition Theorem for $\operatorname{GL}(n, \mathbb{R})$}

\begin{thm}\label{mom}
Let $\operatorname E$ be a linear subspace of $\operatorname{Sym}(n)$ such
that:
$$
[ X, [X, Y]] \in \operatorname E, \qquad for~all\quad X, Y \in \operatorname E,
$$
and let $\operatorname F$ be its orthogonal in $\operatorname{Sym}(n)$:
$$
\operatorname F := \operatorname E^{\perp} = \{ ~X \in \operatorname{Sym}(n)~|~ \Tr XY = 0,~~
\forall Y \in \operatorname E ~\}.
$$
Then for all symmetric positive-definite operator $x$ in
$\operatorname{SPD}(n)$, there exist a unique element $e \in \calE :=\exp
E$ and a unique element $f \in \calF := \exp \operatorname F$ such that $x = e f e.$ More precisely, 
$$
\begin{array}{l}
e := (\pi(x))^{-\frac{1}{2}} \in \calE = \exp \operatorname{E}\\
f := (\pi(x))^{-\frac{1}{2}} x (\pi(x))^{-\frac{1}{2}} \in \calF := \exp \operatorname F,
\end{array}
$$
where $\pi$ denotes the projection onto the totally geodesic submanifold $\exp \operatorname E$.
Moreover the map defined from $\operatorname{SPD}(n)$ to $\exp \operatorname E \times \exp
\operatorname F$ taking $A$ to $(e, f)$ is a homeomorphism.
\end{thm}

\prooft \ref{mom}:\\
Let us introduce the map
$$
\begin{array}{lclc}
\Upsilon: &\exp \operatorname E \times \exp \operatorname F&\longrightarrow &\operatorname{SPD}(n)\\
& (e, f)& \longmapsto & e f e
\end{array}
$$

 Let us show that $\Upsilon$ is one-one. Suppose that $(e_{1},
f_{1})$ and $(e_{2}, f_{2})$ are elements of $\exp \operatorname E \times \exp
\operatorname F$ such that $e_{1} f_{1} e_{1} = e_{2} f_{2} e_{2}$. Consider the
geodesic triangle joining  $e_{1} f_{1} e_{1}$, $e_{1}^{2}$ and
$e_{2}^{2}$. By Theorem \ref{equivgeospace}, $\exp \operatorname E$ is a
geodesic subspace of $\operatorname{SPD}(n)$. Thus the geodesic joining
$e_{1}^{2}$ to $e_{2}^{2}$ lies in $\exp \operatorname E$. On the other hand the
geodesic joining $e_{1} f_{1} e_{1}$ to $e_{1}^{2}$ lies in $e_{1}
\exp \operatorname F \,e_{1}$. Since $\operatorname E$ is perpendicular to $\operatorname F$ at zero, $\exp
\operatorname E$ is perpendicular to $\exp \operatorname F$ at the identity by lemma
\ref{angle}. Since the map taking any $x\in \operatorname{SPD}(n)$ to $e_{1} x\, e_{1}\in \operatorname{SPD}(n)$ is an
isometry, the manifold $e_{1} \exp \operatorname F\, e_{1}$ is perpendicular to the manifold $e_{1}
\exp \operatorname E\, e_{1} = \exp \operatorname E$ at $e_{1}^{2}$. Hence the angle at
$e_{1}^{2}$ of the above geodesic triangle is $90^{\circ}$.
Similarly, the angle at $e_{2}^{2}$ is $90^{\circ}$ since it is
formed by the geodesic joining $e_{2}^{2}$ to $e_{2} f_{2} e_{2} =
e_{1} f_{1} e_{1}$ which lies in $e_{2} \exp F\, e_{2}$ and the
geodesic joining $e_{2}^{2}$ to $e_{1}^{2}$ which lies in $\exp
E$. Denoting by $a$ the length of the side of the geodesic
triangle joining $e_{1}^{2}$ to $e_{2}^{2}$, $b$ the length of the
side joining $e_{1} f_{1} e_{1} $ to $e_{1}^{2}$ and $c$ the
length of the side joining $e_{1} f_{1} e_{1}$ to $e_{2}^{2}$,
ones has $c^{2} \geq b^{2} + a^{2}$ and $b^{2} \geq c^{2} + a^{2}$
by lemma \ref{angle}. This implies that $a = 0$ and $e_{1}^{2} =
e_{2}^{2}$. It follows that $e_{1} = e_{2}$ and $f_{1} = f_{2}$.

 Let us show that $\Upsilon$ is onto. Consider $x$ in
$\operatorname{SPD}(n)$. By Theorem \ref{projection}, the geodesic joining
$x$ to $\pi(x) \in \exp \operatorname E$ is orthogonal to every geodesic
starting at $\pi(x)$ and contained in $\exp \operatorname E$. Denote by $\gamma$
the geodesic satisfying $\gamma(0) = \textrm{Id}\in \operatorname{SPD}(n)$ and $\gamma(1) =
(\pi(x))^{-\frac{1}{2}} x (\pi(x))^{-\frac{1}{2}}$.  Since $y
\mapsto (\pi(x))^{-\frac{1}{2}} y (\pi(x))^{-\frac{1}{2}}$ is an
isometry,  $\gamma$  is orthogonal to every geodesic starting at
the identity and contained in $$(\pi(x))^{-\frac{1}{2}}\cdot  \exp \operatorname E \cdot
(\pi(x))^{-\frac{1}{2}} = \exp \operatorname E.$$ By lemma \ref{angle}, $\gamma$
is tangent to $\operatorname F = \operatorname E^{\perp}$ at the identity and since it is of
the form $t \mapsto \exp t H$ by lemma \ref{formgeo}, we have $H$
in $\operatorname F$. It follows that $\gamma(1) = \exp H =
(\pi(x))^{-\frac{1}{2}} x (\pi(x))^{-\frac{1}{2}}$ is in $\exp \operatorname F$.
Therefore $x = e f e$ with $e := (\pi(x))^{-\frac{1}{2}}$ in $\exp
\operatorname E$ and $f := (\pi(x))^{-\frac{1}{2}} x (\pi(x))^{-\frac{1}{2}}$ in
$\exp \operatorname F$. Consequently $\Upsilon$ is onto.

 The continuity of the map that takes $x$ to $(e, f) \in \exp \operatorname E
\times \exp \operatorname F$ with $x = e f e$, i.e.  $e := (\pi(x))^{-\frac{1}{2}}$  and $f := (\pi(x))^{-\frac{1}{2}} x (\pi(x))^{-\frac{1}{2}}$ follows directly from the
continuity of the projection $\pi$. \cqfdt

\begin{thm}[Mostow's Decomposition]\label{glm}
Let $\operatorname E$ and $\operatorname F$ be as in Theorem \ref{mom}. Then
$\operatorname{GL}(n, \mathbb{R})$ is homeomorphic to the product
$\operatorname{O}(n)\times \exp F\times \exp E$.
\end{thm}

\prooft \ref{glm}:\\
Denote by $\Theta$ the map from $\operatorname{O}(n) \times \exp \operatorname E \times
\exp \operatorname F$ to $\operatorname{GL}(n, \mathbb{R})$ that takes $(k, f, e)$ to the product of matrices $k\cdot f \cdot
e$.

 Let us show that $\Theta $ is one-one. Suppose that $a = k_{1} f_{1}
e_{1} = k_{2} f_{2} e_{2}$ with $(k_{1}, f_{1}, e_{1})$ and
$(k_{2}, f_{2}, e_{2})$ in $\operatorname{O}(n) \times \exp \operatorname E \times
\exp \operatorname F$. We have
$$
a^{T}a = e_{1} f^{2}_{1} e_{1} = e_{2} f_{2}^{2} e_{2}.
$$
Since $f_{1}^{2}$ and $f_{2}^{2}$ are in $\exp \operatorname F$, by Theorem~\ref{mom}, it follows that $e_{1} = e_{2}$ and $f_{1}^{2} =
f_{2}^{2}$. Thus $f_{1} = f_{2}$ and $k_{1} = k_{2}$.

Let us show that $\Theta$ is onto. Consider $x$ in $\operatorname{GL}(n, \mathbb{R})$.
By Theorem \ref{mom},
since $x^{T}x$ is an element of $\operatorname{SPD}(n)$,
there exist $e \in \exp E $ and $f \in \exp F$ such that $x^{T}x =
e f^{2} e$. Let $k$ be $x (f e)^{-1}$. We have:
$$
k^{T}k = \left((f e)^{-1}\right)^T x^{T}x (fe)^{-1} = \left(e^{-1} f^{-1}\right)^T x^{T}x\,e^{-1} f^{-1}  =  f^{-1} e^{-1} \left(e f^{2} e\right)
e^{-1} f^{-1} = \textrm{ id}.
$$
Thus $k$ is in $\operatorname{O}(n)$ and $x = k f e$.

 The continuity of the map that takes $x$ in $\operatorname{GL}(n, \mathbb{R})$ to $(k,
f, e) $ in $\operatorname{O}(n)  \times \exp \operatorname F \times \exp \operatorname E$ follows from
the continuity of the map that takes $x$ to $x^{T}x$ and from
Theorem \ref{mom}. \cqfdt
\\
\\
\demthm \ref{debut}:\\
The Theorem of the Introduction is now a summary of previous considerations. \cqfdt
\\
\\
\demcor \ref{cor_convex}:\\
Consider $N$ a geodesically complete and convex submanifold in $\operatorname{SPD}(n)$. Let $x$ by an arbitrary element in $N$. Then the image of $N$ by the isometry $y \mapsto x^{-\frac{1}{2}}\cdot y \cdot x^{-\frac{1}{2}}$ is a geodesically complete and convex submanifold in $\operatorname{SPD}(n)$ containing $\textrm{Id}$. By Theorem~\ref{debut} it is of the form $\exp \operatorname E$ with $\operatorname E$ satisfying $[X, [X, Y]]\in \operatorname E$.
\cqfd
%
%
\begin{appendix}

\section{The differential of the exponential map of a Lie group $G$}

\subsection{First expression of the differential of the exponential map}

The following Proposition explicits the differential of the
exponential map and is well-known in the theory of
finite-dimensional linear group. For the sake of completeness,
we give below a proof using the canonical connection on the
tangent bundle of a Lie group $G$, which was explained to us by P. Gauduchon. The reader will find the
computation of the differential of the exponential map using
powers series as a consequence of Lemma 1 in \cite{Mostow} (this
computation works as well in the infinite-dimensional setting, see
for instance Proposition 2.5.3 page 116 in \cite{Tum}). See also \cite{God}.

Let us introduce the connection on the tangent bundle $TG$ of
$G$ for which the left-invariant vector fields are parallel.  It is
a flat connection since a trivialization of the tangent bundle is
given by the left-invariant vector fields associated to an arbitrary
basis of the Lie algebra $\g$ of $G$. We give in Proposition~\ref{maurer_connection} the corresponding covariant derivative $\nabla$  in terms of the  (left) Maurer-Cartan $\g$-valued $1$-form $\theta: TG \rightarrow \g$
defined by
$$
\theta_{g}(W) = (L_{g^{-1}})_{*}(W),
$$
where  $g \in G$, $W\in T_gG$, and where $(L_{g^{-1}})_{*}$ denotes the differential of the left multiplication by $g^{-1}$. For a matrix group $G$, $(L_{g^{-1}})_{*}(W) = g^{-1} W$.

\begin{prop}\label{maurer_connection}
Let $W$ be a vector field on $G$, i.e. a section of the tangent bundle $TG$. For a tangent vector $Z$ in
$T_{g}G$, define
\begin{equation}\label{defnabla}
\left(\nabla_{Z}W\right)(g) := (L_{g})_{*}\left({\frac{d}{dt}}_{| t = 0}\theta_{\gamma_Z(t)}\left(W\left(\gamma_Z(t)\right)\right)
\right)\in T_gG
\end{equation}
where $\gamma_Z$ is any smooth curve in $G$ with $\gamma_Z(0) = g$ and ${\frac{d}{dt}}_{| t = 0}\gamma_Z(t) = Z$. 
Then $\nabla$ defines a connection for which left-invariant vector fields are parallel.
\end{prop}

\proofp
Note that $t \mapsto \theta_{\gamma_Z(t)}\left(W\left(\gamma_Z(t)\right)\right)$ is a smooth curve in $\g$, hence its derivative belongs to $\g$. For a left-invariant vector field $W$ this curve is constant, hence its derivative vanishes. In order to verify that formula~\eqref{defnabla} defines indeed a connection, we have to show that 
$$
\nabla_Z(f W) = df(Z) W + f \nabla_Z W,
$$
for any function $f$ on $G$, and any vector fierld $W$. One has
\begin{align*}
\left(\nabla_Z(f W)\right)(g) &= (L_{g})_{*}\left({\frac{d}{dt}}_{| t = 0}\theta_{\gamma_Z(t)}\left(f\left(\gamma_Z(t)\right)W\left(\gamma_Z(t)\right)\right)\right)\\ &  = (L_{g})_{*}\left({\frac{d}{dt}}_{| t = 0}(L_{\gamma_Z(t)^{-1}})_*\left(f\left(\gamma_Z(t)\right)W\left(\gamma_Z(t)\right)\right)\right)\\ 
& = 
 f(g)(L_{g})_{*}\left({\frac{d}{dt}}_{| t = 0}(L_{\gamma_Z(t)^{-1}})_*\left(W\left(\gamma_Z(t)\right)\right)\right)\\
 & + 
 (L_{g})_{*}\left({\frac{d}{dt}}_{| t = 0}(L_{g^{-1}})_*\left(f\left(\gamma_Z(t)\right)W\left(g\right)\right)\right)\\
 & = f(g) (\nabla_Z W)(g) + df(Z) W(g)
 \end{align*}
\cqfd

\begin{prop}\label{log}
For every Lie group $G$, with Lie algebra $\g$, the
differential of the exponential map $~\exp~: \g \rightarrow G$ is
given at $X \in \g$ by~:
\begin{equation}\label{difexp}
\left(d_{X}\exp\right)(Y) = L_{\exp(X)}\left(\frac{1 -
  e^{-\ad(X)}}{\ad(X)}\right)(Y).
\end{equation}
for all $Y$ in $\g$.
\end{prop}

\dem \ref{log}:\\
Let us define the following map~:
$$
\begin{array}{llll}
\Phi~: & \mathbb{R}^2 & \longrightarrow & G\\
       & (t, s) & \longmapsto & \exp\left(t(X + sY)\right).
\end{array}
$$
Consider the push-forward $U$ and $V$ of the vector fields
$\frac{\del}{\del t}$ and $\frac{\del}{\del s}$ on
$\mathbb{R}^2$~:
$$U\left(\Phi(t, s)\right) := \Phi_{*}\left(\frac{\del}{\del t}\right)~~\textrm{and}~~ V\left(\Phi(t,
s)\right) = \Phi_{*}\left(\frac{\del}{\del s}\right).$$ Denote by
$[\cdot\,,\,\cdot]_{\mathfrak{X}}$ the bracket of vector fields.
One has:
\begin{equation}\label{comUV}
[U, V]_{\mathfrak{X}} = \left[\Phi_{*}\left(\frac{\del}{\del
t}\right), \Phi_{*}\left(\frac{\del}{\del
s}\right)\right]_{\mathfrak{X}} = \Phi_{*}\left[\frac{\del}{\del
t}, \frac{\del}{\del s}\right]_{\mathfrak{X}} = 0.
\end{equation}
 Note that
$$
V\left(\Phi(t, s)\right) = \frac{\del \Phi}{\del s}(t, s) =
\left(d\exp_{(tX + stY)}\right)(tY) ~~\textrm{and}~~V\left(\Phi(1,
0)\right) = \left(d_{X}\exp\right)(Y).
$$
The idea of this proof is to explicit the differential equation
satisfied by the $\g$-valued function
$${v}(t):=
\left(L_{\Phi(t, 0)}^{-1}\right)_{*}V\left(\Phi(t, 0)\right) =
\left(L_{\exp(tX)}^{-1}\right)_{*}V\left(\exp(tX)\right).$$
 For this purpose
we will use the connection $\nabla$ on the tangent bundle of
$G$ defined by \eqref{defnabla}.  Note that by the very definition of the exponential
map on a Lie group, $t \mapsto \exp(tX) = \Phi(t, 0)$ is a
geodesic for this connection, hence
\begin{equation}\label{unablau}
\nabla_{U}U = 0
\end{equation}
 along $\Phi(t, 0)$. 
 Let us denote by $T$ and $R$
the torsion and the curvature of $\nabla$. By definition~:
$$
T(U, V) := \nabla_{U}V - \nabla_{V}U - [U, V]_{\mathfrak{X}}
$$
$$\!\!\!\textrm{and}\qquad R_{U, V}U := \nabla_{V}\nabla_{U}U - \nabla_{U}\nabla_{V}U -
\nabla_{[V, U]_{\mathfrak{X}}} U.
$$
By \eqref{comUV}, one has
$$
\nabla_{U}V = \nabla_{V}U  + T(U, V),
$$
hence
$$ \nabla_{U}\left(\nabla_{U}V\right) = \nabla_{U}\nabla_{V}U
+ \nabla_{U}T(U, V).
$$
But the curvature tensor vanishes, hence   \eqref{comUV} and
\eqref{unablau} imply
$$
 \nabla_{U}\nabla_{V}U =  \nabla_{V}\nabla_{U}U -
\nabla_{[V, U]_{\mathfrak{X}}}U  = 0.
$$
Consequently one has
$$
\nabla_{U}\left(\nabla_{U}V\right) = \nabla_{U}T(U, V).
$$
By the expression \eqref{defnabla} of the connection, one has
$$
T(U, V)\left(\Phi(t, 0)\right) = \left( \nabla_{U}V -
\nabla_{V}U\right)(\Phi(t, 0)) = \left(L_{\Phi(t,
0)}\right)_{*}\left(U\cdot\theta(V) - V\cdot\theta(U)\right).
$$
Let us recall that the torsion is a tensor, hence $T(U,
V)\left(\Phi(t, 0)\right)$ does not depend on the extensions of
the vectors $U\left(\Phi(t, 0)\right)$ and $V\left(\Phi(t,
0)\right)$ into vector fields. Using the left-invariant extensions
of these two vectors one see easily that by the very definition of
the bracket in the Lie algebra $\g$ one has
$$
T(U, V)(\Phi(t, 0)) = -\left(L_{\Phi(t,
0)}\right)_{*}\left[\theta(U), \theta(V)\right].
$$
Whence
$$
\nabla_{U}\left(\nabla_{U}V\right) = -\nabla_{U}\left(L_{\Phi(t,
0)}\right)_{*}\left[\theta(U), \theta(V)\right] = \left(L_{\Phi(t,
0)}\right)_{*} \frac{d}{dt}\left[\theta(U), \theta(V)\right].
$$
 Now, along $\Phi(t, 0)$, the vector $\theta(U)$ is the
constant vector $X$, and $\theta(V) = v(t)$. It follows that
$$
\frac{d^{2}v(t)}{dt^2}  = \left(L_{\Phi(t,
0)}^{-1}\right)_{*}\nabla_{U}\left(\nabla_{U}V\right) = -
\left(L_{\Phi(t, 0)}^{-1}\right)_{*}\nabla_{U}\left(L_{\Phi(t,
0)}\right)_{*}\left[X , \theta(V)\right] = - \nabla_{U}[X,
\theta(V)].
$$
This leads to the following differential equation
$$
\frac{d^{2}v(t)}{dt^2}  = - \left[X, \frac{dv}{dt}\right]
$$
with initial conditions $ v(0) = 0 $ and $ \frac{dv}{dt}_{|t=0} =
Y.$ A first integration leads to $$\frac{dv}{dt} =
e^{-t\ad(X)}(Y)$$ and a second to $$v(t) = \left(\frac{1 -
e^{-t\ad(X)}}{\ad(X)}\right)(Y).$$ So the result follows from the
identity $v(1) =
\left(L_{\exp(X)}^{-1}\right)_{*}\left(d_{X}\exp\right)(Y)$.\cqfd

\subsection{Second expression of the differential of the exponential map}

\begin{lem}\label{tx}
For $X$ in $\operatorname{Sym}(n)$, define
\begin{equation}
\begin{array}{llll}
\tau_{X}~: &\operatorname{Sym}(n) &\longrightarrow &
\operatorname{Sym}(n)\\
&  Y & \mapsto &\tau_{X}(Y) := L_{\exp(-\frac{X}{2})}
R_{\exp(-\frac{X}{2})} d_{X}\exp(Y).
\end{array}
\end{equation}
Then 
\begin{equation}\label{secondTx}
\tau_{X} = \frac{\sinh(\ad(X/2))}{\ad(X/2)}  = \sum_{n = 0}^{+\infty}
\frac{(\ad(X)/2)^{2n}}{(2n+1)!}.
\end{equation}
Moreover $\tau_{X} $ is a linear isomorphism of $\operatorname{Sym}(n)$, that depends smoothly on $X\in \operatorname{Sym}(n)$, and whose inverse is 
\begin{equation}
\begin{array}{llll}
\tau_{X}^{-1}:& \operatorname{Sym}(n)   &\rightarrow &\operatorname{Sym}(n) \\
              & Y & \mapsto & \frac{\ad(X/2)}{\sinh (\ad(X/2))}(Y),
\end{array}
\end{equation}
In particular,
\begin{align}\label{second_difexp}
d_{X}\exp(Y)  = R_{\exp(\frac{X}{2})} L_{\exp(\frac{X}{2})}\tau_{X}(Y) = \exp\left(\frac{X}{2}\right)\cdot \left(\frac{\sinh(\ad(X/2))}{\ad(X/2)}(Y)\right)\cdot \exp\left(\frac{X}{2}\right) \\
\end{align}
and
\begin{align}\label{second_difexp2}
(d_{X}\exp)^{-1}(Z) = \tau_X^{-1} L_{\exp(-\frac{X}{2})} R_{\exp(-\frac{X}{2})} (Z) = \frac{\ad(X/2)}{\sinh (\ad(X/2))}\left( \exp\left(-\frac{X}{2}\right)\cdot Z \cdot \exp\left(-\frac{X}{2}\right)\right)
\end{align}

\end{lem}

\proofl \ref{tx}:\\
A direct consequence of formula \eqref{difexp} in
Proposition \ref{log} is that, for all $Y$ in $\operatorname{Sym}(n)$, we have
\begin{align*}
\tau_{X}(Y) & = L_{\exp(-\frac{X}{2})}
R_{\exp(-\frac{X}{2})} d_{X}\exp(Y)  = L_{\exp(-\frac{X}{2})}
R_{\exp(-\frac{X}{2})} L_{\exp(X)}\left(\frac{1 -
  e^{-\ad(X)}}{\ad(X)}\right)(Y) \\ & = L_{\exp(\frac{X}{2})}
R_{\exp(-\frac{X}{2})} \left(\frac{1 -
  e^{-\ad(X)}}{\ad(X)}\right)(Y) = \textrm{Ad}_{\exp(\frac{X}{2})}\left(\left(\frac{1 -
  e^{-\ad(X)}}{\ad(X)}\right)(Y)\right)\\ & = \exp\left(\ad\left(\frac{X}{2}\right)\right)\left(\left(\frac{1 -
  e^{-\ad(X)}}{\ad(X)}\right)(Y)\right) = \left(\frac{e^{\ad(X/2)} -
  e^{-\ad(X/2)}}{\ad(X)}\right)(Y) \\&= \frac{\sinh(\ad(X/2))}{\ad(X/2)}(Y),
\end{align*}
where 
\begin{equation}
\frac{\sinh(\ad(X/2))}{\ad(X/2)} = 
\frac{\exp \left(\ad(X/2)\right) - \exp\left(
-\ad(X/2)\right)}{\ad(X)} = \sum_{n = 0}^{+\infty}
\frac{(\ad(X)/2)^{2n}}{(2n+1)!}.
\end{equation}
Every $X$ in
$\operatorname{Sym}(n)$ is diagonalisable in an orthonormal basis $\{e_i, 1\leq i \leq n\}$ with real eigenvalues $\lambda_{i}$, $1\leq i \leq n$ (counted with multiplicity).  The eigenvalues of $\ad(X)$ acting on $\operatorname{M}(n, \mathbb{R})$ are the real numbers $(\lambda_{i} -
\lambda_{j})$,  $1\leq i, j \leq n$, the correponding eigenvectors being $E_{ij} := e_i\otimes e_j$ (this follows immedatily from the fact that $x$ is similar to a diagonal matrix with eigenvalues $\lambda_i$)). The eigenvalues of $\tau_{X}$ are $$
 \frac{\sinh(\frac{\lambda_{i}-
\lambda_{j}}{2})}{\frac{(\lambda_{i} -
    \lambda_{j})}{2}},$$
    $1\leq  i,j \leq n$, with $E_{ij}$ as eigenvectors (the operator should be interpreted as the identity map when $i=j$). Since for any real number $u\in\mathbb{R}$, $\frac{\sinh u}{u} \geq 1$,   $\tau_{X}$ is injective on $\operatorname{Sym}(n)$. Since $\tau_{X}$ is linear, it is also surjective. Moreover $\tau_{X}$ depends smoothly on $X$ because  $\frac{\sinh u}{u} $ is smooth and $\ad(X)$ depends smoothly on $X$. The inverse of $\tau_{X}$ is the linear map on $\operatorname{Sym}(n)$ with $E_{ij}$ as eigenvectors and corresponding eigenvalues $$
 \frac{\frac{(\lambda_{i} -
    \lambda_{j})}{2}}{\sinh(\frac{\lambda_{i}-
\lambda_{j}}{2})},$$
    $1\leq i, j \leq n$. Therefore it can be written as
$$
\begin{array}{llll}
\tau_{X}^{-1}:& \operatorname{Sym}(n)  &\rightarrow &\operatorname{Sym}(n)\\
              & Y & \mapsto & \frac{\ad(X/2)}{\sinh \ad(X/2)}(Y),
\end{array}
$$
 \cq

\section{Geodesics in locally symmetric homogeneous spaces}
\begin{prop}\label{genial}
Let $M = G/K$ be a locally symmetric homogeneous space of a 
Lie group $G$ with unit element $e$ and Lie algebra $\mathfrak{g} = \mathfrak{k} \oplus
\mathfrak{m}$, where $\mathfrak{k}$ is the Lie algebra of $K$ and $[\mathfrak{m}, \mathfrak{m}] \subset
\mathfrak{k}$. Then, for \textit{any} $G$-invariant Riemannian metric on $M$, the geodesics starting at $o = eK$ are given by the action of one-parameter subgroups of $G$ generated by elements in $\mathfrak{m}$, i.e. are  of the form
$$
\gamma(t) = \left(\exp t\a\right)\cdot o \in M, $$
where $\a$ belongs to $\mathfrak{m}$ (here the dot denotes the action of $G$ on $M$).
\end{prop}

The proof of previous Proposition is based on the fact that for locally symmetric homogeneous manifolds, the homogeneous connection is the Levi-Cevita connection of any $G$-invariant Riemannian metric. The homogeneous connection $\hat{\nabla}$ on
the tangent space of $M$ is defined as follows.
 For every element $\a$ in $\mathfrak{m}_{x}$ and every vector
field $X$ on $M$, one has
 \begin{equation}\label{connexhom}
\hat{\nabla}_{X^{\a}(x)}X = \left(\mathcal{L}_{X^{\a}}X\right)(x)
=  [X^{\a}, X]_{\mathfrak{X}}
 \end{equation}
where $\mathcal{L}$ denotes the Lie derivative and
$[\cdot\,,\,\cdot]_{\mathfrak{X}}$ the bracket of vector fields.
 For $\a$ in
$\mathfrak{m}_{x}$ and $\b$ in $\g$, one has~:
$$
\hat{\nabla}_{X^{\a}(x)}X^{\b} =  [X^{\a},
X^{\b}]_{\mathfrak{X}}(x) = - X^{[\a, \b]}(x).
$$
The torsion of the connection $\hat{\nabla}$ is given by
$$
T^{\hat{\nabla}}(X^{\a}, X^{\b}) = \hat{\nabla}_{X^{\a}}X^{\b} -
\hat{\nabla}_{X^{\b}}X^{\a} - [X^{\a}, X^{\b}]_{\mathfrak{X}} = -
X^{[\a, \b]}.
$$
It follows that for a locally symmetric homogeneous space, the
homogeneous connection is torsion free since for $\a$ and $\b$ in
$T_{x}M = \mathfrak{m}_{x}$, $[\a, \b]$ belongs to the isotropy
$\mathfrak{k}_{x}$ thus $X^{[\a, \b]}$ vanishes. On the other
hand, it follows from definition \eqref{connexhom} that the
covariant derivative of any tensor field $\Phi$ along $Y \in
T_{x}M$ is the Lie derivative of $\Phi$ along the vector field
$X^{\a}$ where $\a \in \mathfrak{m}_{x}$ is such that $Y =
X^{\a}(x)$. Thus the homogeneous connection preserves every
$G$-invariant Riemannian metric. Consequently $\hat{\nabla}$ is
the Levi-Civita connection of every $G$-invariant Riemannian
metric on the locally symmetric space $M = G/H$.

\vspace{0.5cm}
\dem \ref{genial}:\\
Every element $\a$ in $\g$ generates a vector field $X^{\a}$ on
the homogeneous space $M = G/K$.
For every $x = g\cdot o$, $g \in G$, the Lie algebra $\g$ splits
into $\g = \mathfrak{k}_{x} \oplus \mathfrak{m}_{x}$, where
$\mathfrak{k}_{x} := \textrm{Ad}(g)(\mathfrak{k})$ is the Lie
algebra of the isotropy group at $x$ and where $\mathfrak{m}_{x}
:= \textrm{Ad}(g)(\mathfrak{m})$ can be identified with the
tangent space $T_{x}M$ of $M$ at $x$ by the application $\a
\mapsto X^{\a}(x)$. To see that for $\a \in \mathfrak{m}$, the
curve
$$ \gamma(t) = \left(\exp
t\a\right)\cdot o, \qquad \a \in \mathfrak{m}.
$$
is a geodesic, note that the equality $\a = \textrm{Ad}(\exp
t\a)(\a)$ implies that $\a$ belongs to the space
$\mathfrak{m}_{\gamma(t)}$ for all $t$. Hence from
$\dot{\gamma}(t) = X^{\a}(\gamma(t))$ it follows that
$\hat{\nabla}_{\dot{\gamma}(t)}\dot{\gamma}(t) =
\mathcal{L}_{X^{\a}}X^{\a}(\gamma(t)) = 0$. In other words
$\gamma$ is a geodesic of $M$. \cqfd

\medskip

\begin{rem}[Curvature] 
For a locally symmetric homogeneous space, the curvature tensor $R$ of $\hat{\nabla}$ defined by 
$$
R_{X, Y}Z := \hat{\nabla}_{X}\hat{\nabla}_{Y}Z - \hat{\nabla}_{Y} \hat{\nabla}_{X}Z -
\hat{\nabla}_{[X, Y]} Z.
$$
has the following simple formula
\begin{equation}\label{curvature_symmetric}
R_{X, Y}Z = [[X, Y], Z].
\end{equation}
Indeed, for $\a, \b, \mathfrak{c} \in \mathfrak{m}_{x}$, one has
\begin{align*}
R_{X^\a, X^\b} X^{\mathfrak{c}} & = \hat{\nabla}_{X^\a}( [X^{\b}, X^{\mathfrak{c}}]_{\mathfrak{X}}) 
- \hat{\nabla}_{X^\b}( [X^{\a}, X^{\mathfrak{c}}]_{\mathfrak{X}}) 
- \hat{\nabla}_{[X^{\a}, X^{\b}]_\mathfrak{X}} X^{\mathfrak{c}} 
\\
&
=[X^{\a}, [X^{\b}, X^{\mathfrak{c}}]_{\mathfrak{X}}]_{\mathfrak{X}}
- [X^{\b}, [X^{\a}, X^{\mathfrak{c}}]_{\mathfrak{X}}]_{\mathfrak{X}}
+ \hat{\nabla}_{X^{[{\a}, {\b}]}} X^{\mathfrak{c}} \\
&
= [[X^{\a}, X^{\b}]_{\mathfrak{X}}, X^{\mathfrak{c}}]_{\mathfrak{X}},
\end{align*}
where in the last equality we have used $X^{[{\a}, {\b}]} = 0$ for $\a, \b\in \mathfrak{m}_{x}$, since $[\mathfrak{m}_{x}, \mathfrak{m}_{x}]\subset \mathfrak{k}_x$ and $\mathfrak{k}_x$ acts trivially on $T_xM$.
\end{rem}

\end{appendix}

\vspace{0.5 cm}
\begin{center}
\begin{tabular}{ll}
Alice Barbara TUMPACH & Gabriel Larotonda \\ 
Institut CNRS Pauli &  Departamento de Matem\'atica, FCEyN, UBA\\
Oskar-Morgenstern-Platz 1 & Avenida Cantilo S/N, Ciudad Universitaria\\
1090 Vienna, Austria & 1428 Buenos Aires, Argentina \\ {\it E-mail address~:} {\tt barbara.tumpach@math.cnrs.fr}  & {\tt glaroton@dm.uba.ar}
\end{tabular}
\end{center}
\end{document}